\documentstyle[amstex,amsfonts,amssymb,latexsym,graphicx,bbm,texdraw,epsfig]
{jsc}

\input txdtools.tex

\newcommand{\im}{\text{Im}\,}
\newcommand{\Ker}{\text{Ker}\,}
\newcommand{\Mat}{\text{Mat}\,}
\newcommand{\Hom}{\text{Hom}}
\newcommand{\Ext}{\text{Ext}}
\newcommand{\Ass}{\text{Ass}\,}
\newcommand{\minAss}{\text{minAss}\,}
\newcommand{\Ann}{\text{Ann}\,}
\newcommand{\Spec}{\text{Spec}\,}
\newcommand{\Lc}{\text{lc}\,}
\newcommand{\codim}{\text{codim}\,}
\newcommand{\lm}{\text{lm}\,}
\newcommand{\NF}{\text{NF}\,}
\newcommand{\ann}{\text{ann}\,}
\newcommand{\Loc}{\text{Loc}\,}
\newcommand{\syz}{\text{syz}\,}
\newcommand{\mult}{\text{mult}\,}
\newcommand{\N}{{\mathbb N}}
\newcommand{\R}{{\mathbb R}}
\newcommand{\C}{{\mathbb C}}
\newcommand{\Q}{{\mathbb Q}}
\renewcommand{\P}{{\mathbb P}}
\newcommand{\kg}{{\mathcal G}}
\newcommand{\ko}{{\mathcal O}}
\newcommand{\fm}{\mathfrak{m}}
\newcommand{\lra}{\longrightarrow}
\newcommand{\lc}{lc}
\newcommand{\lcm}{lcm}

\def\cavec (#1 #2)(#3 #4)(#5 #6){
   \clvec (#1 #2)(#3 #4)(#5 #6)
   \cossin (#3 #4)(#5 #6)\cosa\sina
   \rmove (0 0)
   \bsegment
      \drawdim in \setsegscale 0.05
      \move ({-\cosa} -\sina) \lpatt()\avec (0 0)
   \esegment}

\title[Computer Algebra and Algebraic Geometry]{{Computer Algebra and
    Algebraic}\\ Geometry --- 
Achievements and Perspectives*}
\author[Gert-Martin Greuel]{{GERT-MARTIN GREUEL} \\
Fachbereich Mathematik\\
Universit\"at Kaiserslautern,
Erwin-Schr\"odinger-Stra{\ss}e\\
D -- 67663 Kaiserslautern\\
email (greuel@@mathematik.uni-kl.de)
}

\date{\today}
\begin{document}

\maketitle

\begin{flushright}
De computer is niet de steen\\
 maar de slijpsteen der wijzen.\\
(The computer is not the philosopher's stone\\
 but the philosopher's
whetstone.)\\
Hugo Battus, Rekenen op taal (1989)
\end{flushright}

\tableofcontents

\vfill

\noindent * {\footnotesize Extended version of invited talk, delivered at the
            ISSAC'98 conference at Rostock,\\ 
\hspace*{0.4cm} 13th -- 15th August, 1998.}

\section{Preface}

\noindent In this survey I should like to introduce some concepts of algebraic geometry
and try to demonstrate the fruitful interaction between algebraic geometry and
computer algebra and, more generally, between mathematics and computer
science.  One of the aims of this article is to show, by means of examples,
the usefulness of computer algebra to mathematical research.

Computer algebra itself is a highly diversified discipline with applications
to various areas of mathematics; we find many of these in numerous research
papers, in proceedings or in textbooks (cf.\ \ccite{BWi,CCS,MGH,ISSAC}).  Here
we concentrate mainly on Gr\"obner bases and leave aside many other topics of
computer algebra (cf.\ \ccite{DST,GG,GKW}).  In particular, we do not mention
(multivariate) polynomial factorisation, another major and important tool in
computational algebraic geometry.    Gr\"obner bases were introduced
originally by Buchberger as a computational tool for testing solvability of a
system of polynomial equations, to count the number of solutions (with
multiplicities) if this number is finite and, more algebraically, to compute
in the quotient ring modulo the given polynomials.  Since then, Gr\"obner
bases have become \textit{the}\/ major computational tool, not only in
algebraic geometry.

The importance of Gr\"obner bases for mathematical research in algebraic
geometry is obvious and their use needs, nowadays, hardly any justification.
Indeed, chapters on Gr\"obner bases and \hbox{Buchberger's} algorithm
\cite{Bu1} have been incorporated in many new textbooks on algebraic geometry
such as the books of \ccite{CLO1,CLO2} or the recent books
of \ccite{Ei2,Va2}, not to mention textbooks which are devoted
exclusively to Gr\"obner bases, as \ccite{AL,BWe,Fr}.

Computational methods become increasingly important in pure mathematics and
the above mentioned books have the effect that Gr\"obner bases and their
applications become a standard part of university courses on algebraic
geometry and commutative algebra.  One of the reasons is that these methods,
together with very efficient computers, allow the treatment of non--trivial
examples and, moreover, are applicable to non--mathematical, industrial,
technical or economical problems.  Another reason is that there is a belief
that algorithms can contribute to a deeper understanding of a problem.  The
human idea of ``understanding'' is clearly part of the historical, cultural
and technical status of the society and nowadays understanding in mathematics
requires more and more algorithmic treatment and computational mastering.

On the other hand, it is also obvious that many of the recent deepest
achievements in algebraic and arithmetic geometry, such as string
theory and mirror symmetry (coming from physics) or Wiles' proof of
Fermat's last theorem, just to mention a few,  were neither inspired by
nor used computer algebra at all. I just mention this in order to
stress that no computer algebra system can ever replace, in any
significant way, mathematical thinking. 

Generally speaking, algorithmic treatment and computational mastering marks
not the beginning but the end of a development and already requires an
advanced theoretical understanding. In many cases an algorithm is, however,
much more than just a careful analysis of known results, it is really a new
level of understanding, and an efficient implementation is, in  addition,
usually a highly nontrivial task.  Furthermore, having a computer algebra
system which has such algorithms implemented and which is easy to use, it then
becomes a powerful tool for the working mathematician, like a calculator for
the engineer.

In this connection I should like to stress that having Buchberger's algorithm
for computing Gr\"obner bases of an ideal is, although indispensable, not much more than having +,-,*,/ on a
calculator. Nowadays there exist efficient implementations of very involved
and sophisticated algorithms (most of them use Gr\"obner bases in an essential
way) allowing the computation of such things as

\begin{itemize}
\item Hilbert polynomials of graded ideals and modules, 
\item free resolutions of finitely generated modules,
\item Ext, Tor and cohomology groups,
\item infinitesimal deformations and obstructions of varieties and
singularities, 
\item versal deformations of varieties and singularities,
\item primary decomposition of ideals,
\item normalisation of affine rings,
\item invariant rings of finite and reductive groups,
\item Puiseux expansion of plane curve singularities, 
\end{itemize}

\noindent not to mention the standard operations like ideal and radical
membership, ideal intersection, ideal quotient and elimination of variables.
All the above--mentioned algorithms are implemented in SINGULAR \cite{GPS},
some of them also in CoCoA \cite{CNR} and Macaulay, resp.\ Macaulay2 (Bayer
and Stillmann, resp.\ Grayson and Stillmann), to mention computer algebra
systems which are designed for use in algebraic geometry and commutative
algebra.   Even general purpose and commercial systems such as Mathematica,
Maple, MuPad etc.\ offer Gr\"obner bases and, based on this, libraries
treating special problems in algebra and geometry.

It is well--acknowledged that Gr\"obner bases and \hbox{Buchberger's}
algorithm are responsible for the possibility to compute the above objects in
affine resp.\ projective geometry, that is, for non--graded resp.\ graded
ideals and modules over polynomial rings. It is, however, much less known that
standard bases (``Gr\"obner bases`` for not necessarily well--orderings) can
compute the above objects over the localisation of polynomial rings. This is
basically due to Mora's modification of Buchberger's algorithm (\ccite{Mo})
which has been modified and extended to arbitrary (mixed) monomial orderings
in SINGULAR since 1990 and was published in \ccite{Getal} and in \ccite{GP1}.
We include a brief description in Section 5.

I shall explain how non--well--orderings are intrinsically associated with a
ring which may be, for example, a local ring or a tensor product of a local
and a polynomial ring. These ``mixed rings'' are by no means exotic but are
necessary for certain algorithms which use tag--variables which have to be
eliminated later. The extension of \hbox{Buchberger's} algorithm to
non--well--orderings has important applications to problems in local algebraic
geometry and singularity theory, such as the computation of

\begin{itemize}
\item local multiplicities,
\item Milnor and Tjurina numbers, 
\item syzygies and Hilbert--Samuel functions for local rings 
\end{itemize}
but also to more advanced
algorithms such as 
\newpage

\begin{itemize}
\item classification of singularities,
\item semi--universal deformation of singularities,
\item computation of moduli spaces, 
\item monodromy of the Gau{\ss}--Manin connection. 
\end{itemize}

Moreover, I demonstrate, by means of examples, how some of the
above algorithms were used to support mathematical research in a
non--trivial manner. These examples belong to the main methods of
applying computer algebra successfully:
 
\begin{itemize}
\item producing counter examples or giving support to conjectures,
\item providing evidence and prompting proofs for new theorems,
\item constructing interesting explicit examples.
\end{itemize}

The mathematical problems I present were, to a large extent, responsible for
the development of SINGULAR, its functionality and speed.

Finally, I point out some open problems in mathematics and
non--mathematical applications which are a challenge to computer algebra and
where either the knowledge of an algorithm or an efficient implementation is
highly desirable.
\smallskip

\textbf{Acknowledgement}: The author was partially supported by the DFG
Schwerpunkt ``Effiziente Algorithmen f\"ur diskrete Probleme''.  Special thanks
to C.~Lossen and, in particular, to T.~Keilen for preparing the pictures.
Finally, I should like to thank the referees for useful comments.

\section{Introduction by pictures}

The basic problem of algebraic geometry is to understand the set of points $x =
(x_1, \dots, x_n) \in K^n$ satisfying a system of equations
\[
\begin{array}{ccc}
f_1(x_1, \dots, x_n) & = & 0,\\
\vdots & & \\
f_k(x_1, \dots, x_n) & = & 0,
\end{array}
\]
where $K$ is a field and $f_1, \dots, f_k$ are elements of the polynomial ring
$K[x] = K[x_1, \dots, x_n]$.

The solution set of $f_1 = 0, \dots, f_k= 0$ is called the algebraic set, or
\textbf{algebraic variety} of $f_1, \dots, f_k$ and is denoted by
\[
V = V(f_1, \dots, f_k).
\]

It is easy to see, and important to know, that $V$ depends only on the ideal
\[
I = \langle f_1, \dots, f_k\rangle = \{f \in K[x]\mid f = \sum^k_{i=1} a_i f_i,\;
a_i \in K[x]\}
\]
generated by $f_1, \dots, f_k$ in $K[x]$, that is $V = V(I) = \{x \in K^n\mid
f(x) = 0\; \forall\; f \in I\}$.

   \small
   \unitlength1cm
   \begin{picture}(13.5,26)
     \put(0,19.5){\includegraphics[clip,width=6.5cm]{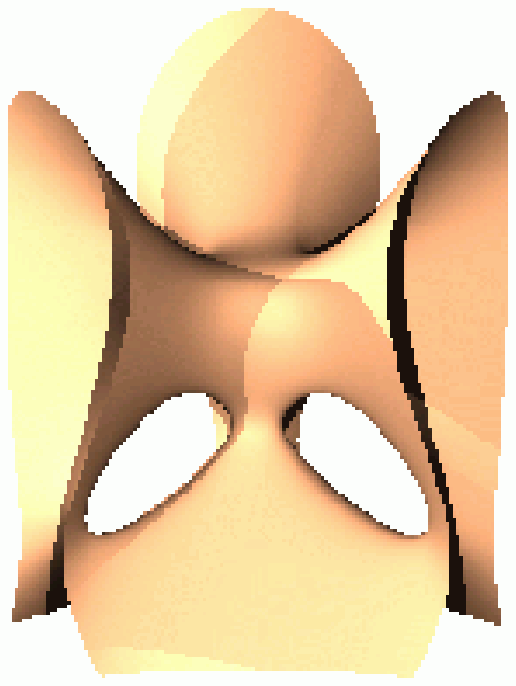}}
     \put(0,19){
       \begin{minipage}[t]{6.5cm}
         \begin{center}
           \textbf{The Clebsch Cubic}
         \end{center}
         This is the unique cubic surface which has ${\mathbb S}_5$,
         the symmetric group of 5 letters, as symmetry group. It is
         named after its discoverer Alfred Clebsch and has the affine equation
         \footnotesize
         \begin{gather*}
           81(x^3+y^3+z^3)\\
           -189(x^2y+x^2z+xy^2+xz^2+y^2z+yz^2)\\
           +54xyz+126(xy+xz+yz)\\
           -9(x^2+y^2+z^2)-9(x+y+z)+1=0. 
         \end{gather*}
       \end{minipage}
       }
     \put(7,19.5){\includegraphics[clip,width=6cm]{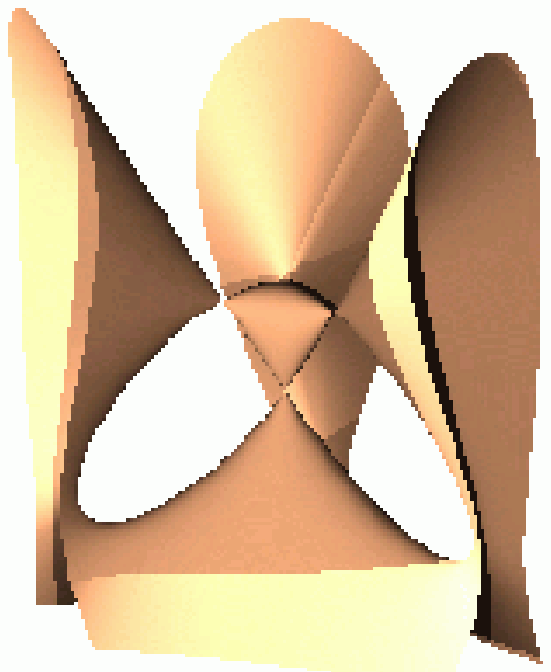}}
     \put(7,19){
       \begin{minipage}[t]{6cm}
         \begin{center}
           \textbf{The Cayley Cubic}
         \end{center}
         There is a unique cubic surface which has four ordinary
         double points, usually called the Cayley cubic after its
         discoverer, Arthur Cayley. It is a degeneration of the Clebsch cubic,
         has ${\mathbb S}_4$
         as symmetry group, and the projective equation is
         \footnotesize
         \begin{displaymath}
           z_0z_1z_2+z_0z_1z_3+z_0z_2z_3+z_1z_2z_3=0.
         \end{displaymath}
       \end{minipage}
       }

     \put(0,10){\includegraphics[clip,width=5.2cm]{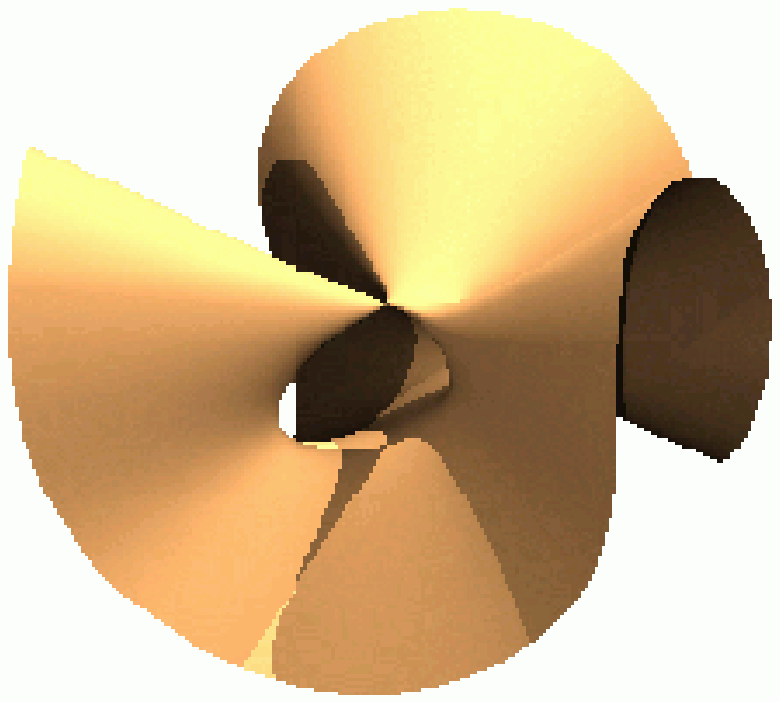}}
     \put(0,8.5){
       \begin{minipage}[t]{6cm}
         \begin{center}
           \textbf{A Cubic with a $D_4$-Singularity}
         \end{center}
         Degenerating the Cayley cubic we receive a
         $D_4$-singularity. The affine equation is
         \footnotesize
         \begin{displaymath}
           x(x^2-y^2)+z^2(1+z)+\tfrac{2}{5}xy+\tfrac{2}{5}yz = 0.
         \end{displaymath}
       \end{minipage}
       }
     \put(7,9.5){\includegraphics[clip,width=6cm]{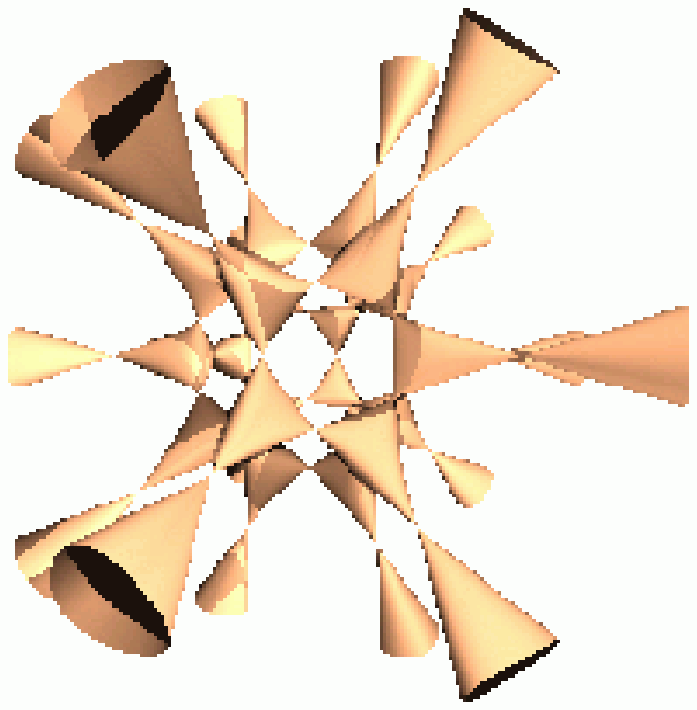}}
     \put(7,8.5){
       \begin{minipage}[t]{6cm}
         \begin{center}
           \textbf{The Barth Sextic}
         \end{center}
         The equation for this sextic was found by Wolf Barth. It has
         65 ordinary double points, the maximal possible number for a
         sextic. Its affine equation is (with $c=\tfrac{1+\sqrt{5}}{2}$)
         \footnotesize
         \begin{gather*}
           (8c+4)x^2y^2z^2 -c^4(x^4y^2+y^4z^2+x^2z^4)\\
           +c^2(x^2y^4+y^2z^4+x^4z^2)\\
           -\frac{2c+1}{4}(x^2+y^2+z^2-1)^2=0. 
         \end{gather*}
       \end{minipage}
       }
   \end{picture}

\newpage

   \small
   \unitlength1cm
   \begin{picture}(13,17)
       \put(0,17){\includegraphics[clip,width=5.5cm,angle=180]{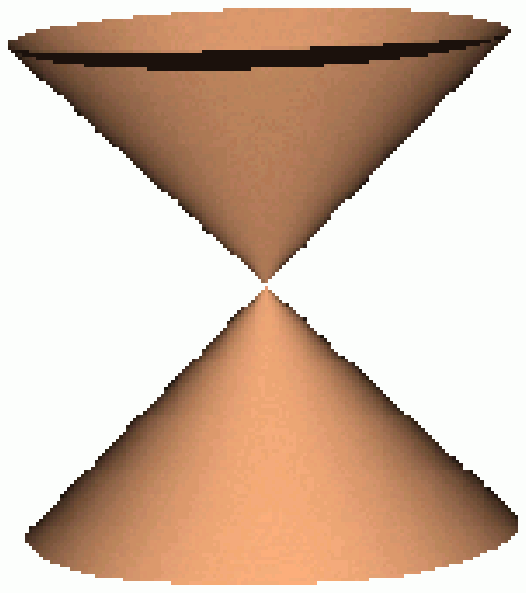}}
      \put(0,10.5){ 
       \begin{minipage}[t]{6cm}
         \begin{center}
           \textbf{An Ordinary Node}
         \end{center}
         An ordinary node is the most simple singularity. It has the
         local equation
         \footnotesize
         \begin{displaymath}
           x^2+y^2-z^2=0.
         \end{displaymath}
       \end{minipage}
       }

     \put(7,11){\includegraphics[clip,width=6cm]{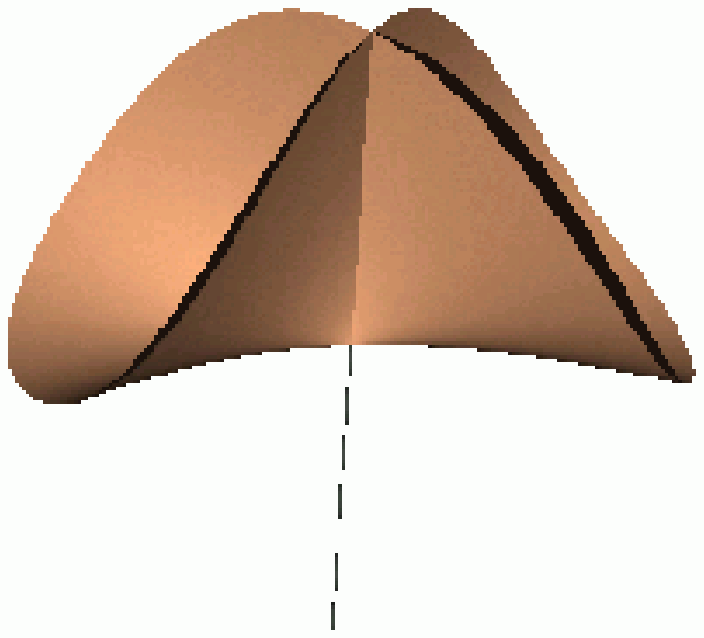}}
       \put(7,11)
       {\begin{minipage}[t]{6cm}
         \begin{center}
           \textbf{Whitney's Umbrella}
         \end{center}
         The Whitney umbrella is named after Hassler Whitney who
         studied it in connection with the stratification of analytic
         spaces. It has the local equation
         \begin{center}$y^2-zx^2=0$.\end{center}
                \end{minipage}
       }

\put(0,0){
\begin{minipage}[b]{6cm}
\begin{center}
 {\includegraphics[trim=20 20 20 20, clip]{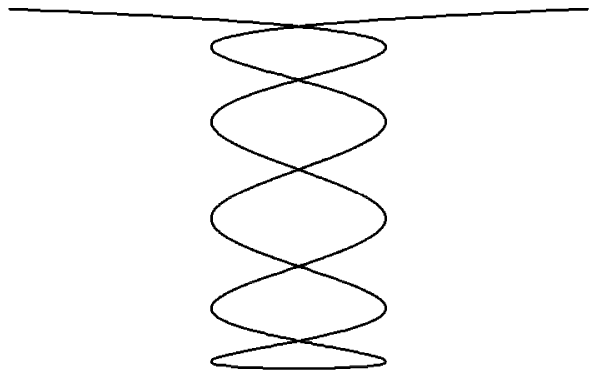}}\\
         A $5$-nodal plane curve of degree $11$\\
 with equation $-16x^2+1048576y^{11}-720896y^9$\\
+$180224y^7
-19712y^5+880y^3-11y+\frac{1}{2}$,\\
 a deformation of 
$A_{10}\; : \;y^{11}-x^2=0$.
\end{center}
       \end{minipage}
       }
\medskip

      \put(7,0){
       \begin{minipage}[b]{6cm}
\begin{center}
{\includegraphics[clip,height=5.5cm,width=5.5cm]{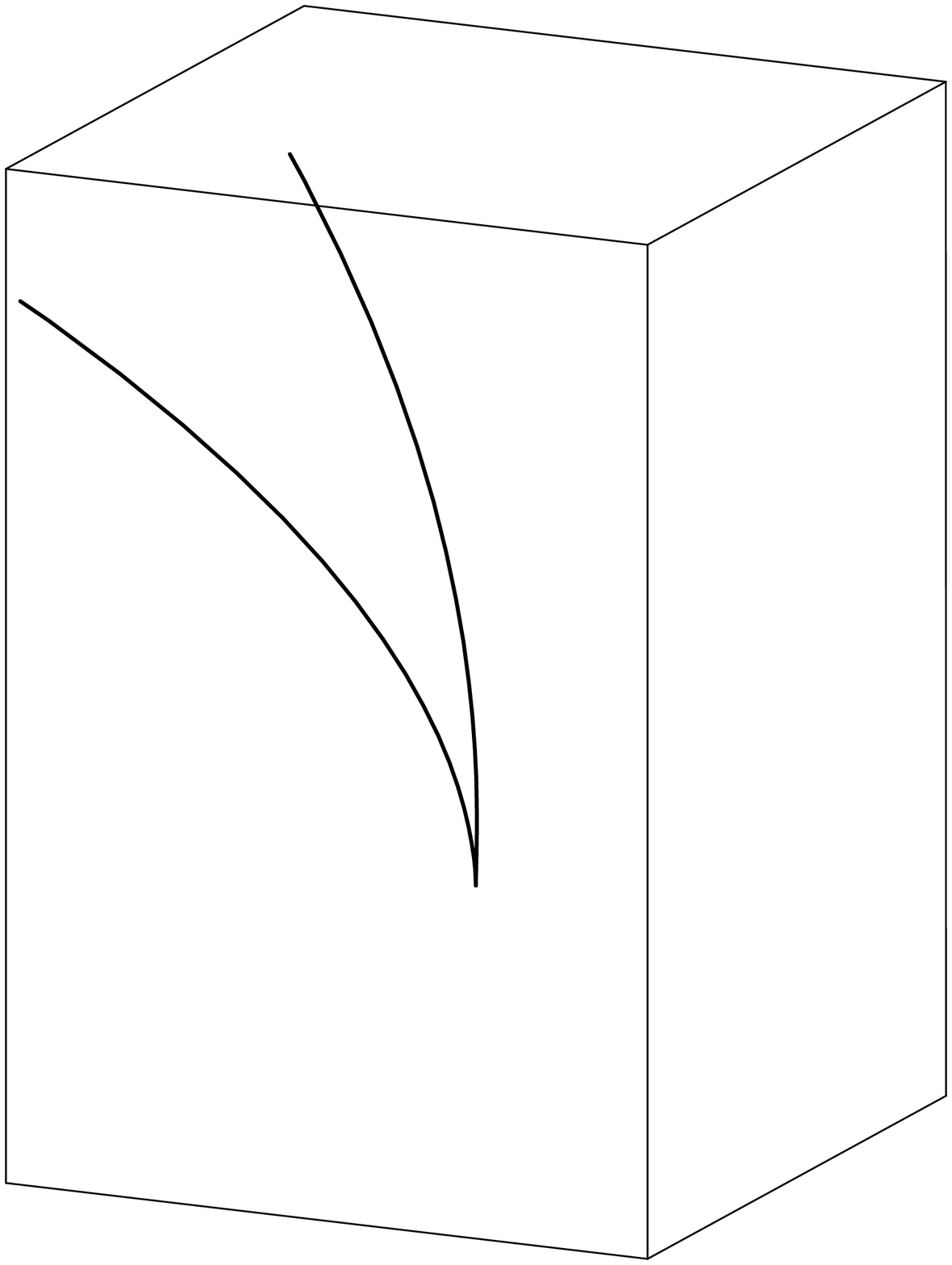}}\\
         This space curve is given parametrically by
        $x=t^4,\; y=t^3,\; z=t^2$,
                  or implicitly by $x-z^2=y^2-z^3=0$.
\end{center}
       \end{minipage}
       }

\end{picture}
\newpage

Of course, if for some polynomial $f \in K[x]$, $f^d|_V = 0$, then $f|_V = 0$
and hence, $V = V(I)$ depends only on the \textbf{radical} of $I$,
\[
\sqrt{I} = \{f \in K[x] \mid f^d \in I, \text{ for some } d\}.
\]
The biggest ideal determined by $V$ is
\[
I(V) = \{f \in K[x] \mid f(x) = 0\; \forall\; x \in V\},
\]
and we have $I \subset \sqrt{I} \subset I (V)$ and $V\bigl(I(V)\bigr) =
V(\sqrt{I}) = V(I) = V$.

The important \textbf{Hilbert Nullstellensatz} states that, for $K$ an algebraically
closed field, we have for any variety $V \subset K^n$ and any ideal $J \subset
K[x]$,
\[
V = V(J) \Rightarrow I(V) = \sqrt{J}
\]
(the converse implication being trivial).  That is, we can recover the ideal
$J$, up to radical, just from its zero set and, therefore, for fields like
$\C$ (but, unfortunately, not for $\R$) geometry and algebra are ``almost
equal''.  But almost equal is not equal and we shall have occasion to see that
the difference between $I$ and $\sqrt{I}$ has very visible geometric
consequences.

Many of the problems in algebra, in particular, computer algebra, have a
geometric origin.  Therefore, I choose an introduction by means of some
pictures of algebraic varieties, some of them being used to illustrate
subsequent problems.

The above pictures were not only chosen to illustrate the beauty of algebraic
geometric objects but also because these varieties have had some prominent
influence on the development of algebraic geometry and singularity theory.

The Clebsch cubic itself has been the object of numerous investigations in
global algebraic geometry, the Cayley and the $D_4$--cubic also, but,
moreover, since the $D_4$--cubic deforms, via the Cayley cubic, to the Clebsch
cubic, these first three pictures illustrate deformation theory, an important
branch of (computational) algebraic geometry.

The ordinary node, also called $A_1$--singularity (shown as a surface
singularity) is the most simple singularity in any dimension.  The Barth
sextic illustrates a basic but very difficult and still (in general) unsolved
problem: to determine the maximum possible number of singularities on a
projective variety of given degree.  In Section 7.3 we report on recent
progress on this question for plane curves.

Whitney's umbrella was, at the beginning of stratification theory, an
important example for the two Whitney conditions.  We use the umbrella in
Section 4.2 to illustrate that the algebraic concept of normalisation may even
lead to a parametrisation of a singular variety, an ultimate goal in many
contexts, especially for graphical representations.  In general, however, such
a parametrisation is not possible, even not locally, if the variety has
dimension bigger than one.  For curve singularities, on the other hand, the
normalisation is always a parametrisation.  Indeed, computing the
normalisation of the ideal given by the implicit equations for the space curve
in the last picture, we obtain the given parametrisation.  Conversely, the
equations are derived from the parametrisation by eliminating $t$, where
elimination of variables is perhaps the most important basic application of
Gr\"obner bases.

Finally, the 5--nodal plane curve illustrates the global existence problem
described in Section 7.2.  Moreover, these kind of deformations with the
maximal number of nodes play also a prominent role in the local theory of
singularities.  For instance, from this real picture we can read off the
intersection form and, hence, the monodromy of the singularity $A_{10}$ by a
beautiful theory of A'Campo and Gusein-Zade.  We shall present a completely
different, algebraic algorithm to compute the monodromy in Section 6.3.

For more than a hundred years, the connection between algebra and geometry has
turned out to be very fruitful and both merged to one of the leading areas in
mathematics: algebraic geometry. The relationship between both disciplines can
be characterised by saying that algebra provides rigour while geometry
provides intuition.

In this connection, I place computer algebra on top of rigour, but I should
like to stress its limited value if it is used without intuition.

\section{Some problems in algebraic geometry}

In this section I shall formulate some of the basic questions and problems
arising in algebraic geometry and provide ingredients for certain algorithms.
I shall restrict myself to those algorithms where I am somehow
familiar with their implementations and which have turned out to be useful in
practical applications.

Let me first recall the most basic but also most important applications of
Gr\"obner bases to algebraic constructions (called ``Gr\"obner basics'' by
Sturmfels).  Since these can be found in more or less any textbook dealing
with Gr\"obner bases, I just mention them:

\begin{itemize}
\item Ideal (resp.\ module) membership problem
\item Intersection with subrings (elimination of variables)
\item Intersection of ideals (resp.\ submodules)
\item Zariski closure of the image of a map
\item Solvability of polynomial equations
\item Solving polynomial equations
\item Radical membership
\item Quotient of ideals
\item Saturation of ideals
\item Kernel of a module homomorphism
\item Kernel of a ring homomorphism
\item Algebraic relations between polynomials
\item Hilbert polynomial of graded ideals and modules
\end{itemize}

The next questions and problems lead to algorithms which are slightly
more (some of them much more) involved.   They are, nevertheless, still
very basic and quite natural.   I should like to illustrate them by
means of four simple examples, shown in the pictures of this section, referred
to as Example 1) -- 4):

Assume we are given an ideal $I \subset K[x_1, \dots, x_n]$ by a set of
generators $f_1, \dots, f_k \in K[x]$. Consider the following questions and
problems:

\begin{enumerate}
\item \textit{Is $V(I)$ irreducible or may it be decomposed into several
    algebraic varieties?   If so, find its irreducible components.
    Algebraically this means to compute a primary decomposition of $I$ or of
    $\sqrt{I}$, the latter means to compute the associated prime ideals of
    $I$.} 
  
  Example 1) is irreducible, Example 2) has two components (one of
  dimension 2 and one of dimension 1), Example 3) has three (one--dimensional)
  and Example 4) has nine (zero--dimensional) components.

\item \textit{Is $I$ a radical ideal (that is, $I = \sqrt{I})$?   If not,
    compute its radical $\sqrt{I}$.}
  
  In Examples 1) -- 3) $I$ is radical while in Example 4) 
  $\sqrt{I} = \langle y^3 -y, x^3-x\rangle$, which is much simpler than $I$.
  In this example the central point corresponds to $V(\langle x,y\rangle^2)$
  which is a \textbf{fat point}, that is, it is a solution of $I$ of
  multiplicity ($= \dim_K K[x,y]/\langle x,y\rangle^2$) bigger than 1 (equal
  to 3).  All other points have multiplicity 1, hence the total number of
  solutions (counted with multiplicity) is 11.  This is a typical example of
  the kind Buchberger (resp.\ Gr\"obner) had in mind at the time of writing
  his thesis.

\item \textit{A natural question to ask is ``how independent are the
    generators $f_1, \dots, f_k$ of $I$?'', that is, we ask for all relations
\[
(r_1, \dots, r_k) \in K[x]^k, \text{ such that } \sum r_i f_i = 0.
\]
 }

These relations form a submodule of $K[x]^k$, which is called the
\textbf{syzygy module} of $f_1, \dots, f_k$ and is denoted by $\syz(I)$.  It
is the kernel of the $K[x]$--linear map
\[
K[x]^k \lra K[x];\;\; (r_1, \dots, r_k) \longmapsto \sum r_i f_i.
\]

\item \textit{More generally, we may ask for generators of the kernel of a
    $K[x]$--linear map $K[x]^r \lra K[x]^s$ or, in other words, for solutions
    of a system of linear equations over $K[x]$.}

A direct geometric interpretation of syzygies is not so clear, but there are
instances where properties of syzygies have important geometric consequences
cf.\ \ccite{Sch3}.

In Example 1) we have $\syz(I) = 0$, in Example 2), $\syz(I) =
\langle(-y,x)\rangle \subset K[x]^2$, in Example 3), $\syz(I) =
\langle(-z,y,0),(-z,0,x)\rangle \subset K[x]^3$ and in Example 4),
$\syz(I)\subset K[x]^4$ is generated by $(x,-y,0,0),
(0,0,x,-y),(0,x^2-1,-y^2+1,0)$. 

\begin{center}
  \textbf{Four examples}
\end{center}

\unitlength1cm
\begin{picture}(13,5.5)
 \put(-1,0.5){\epsfig{file=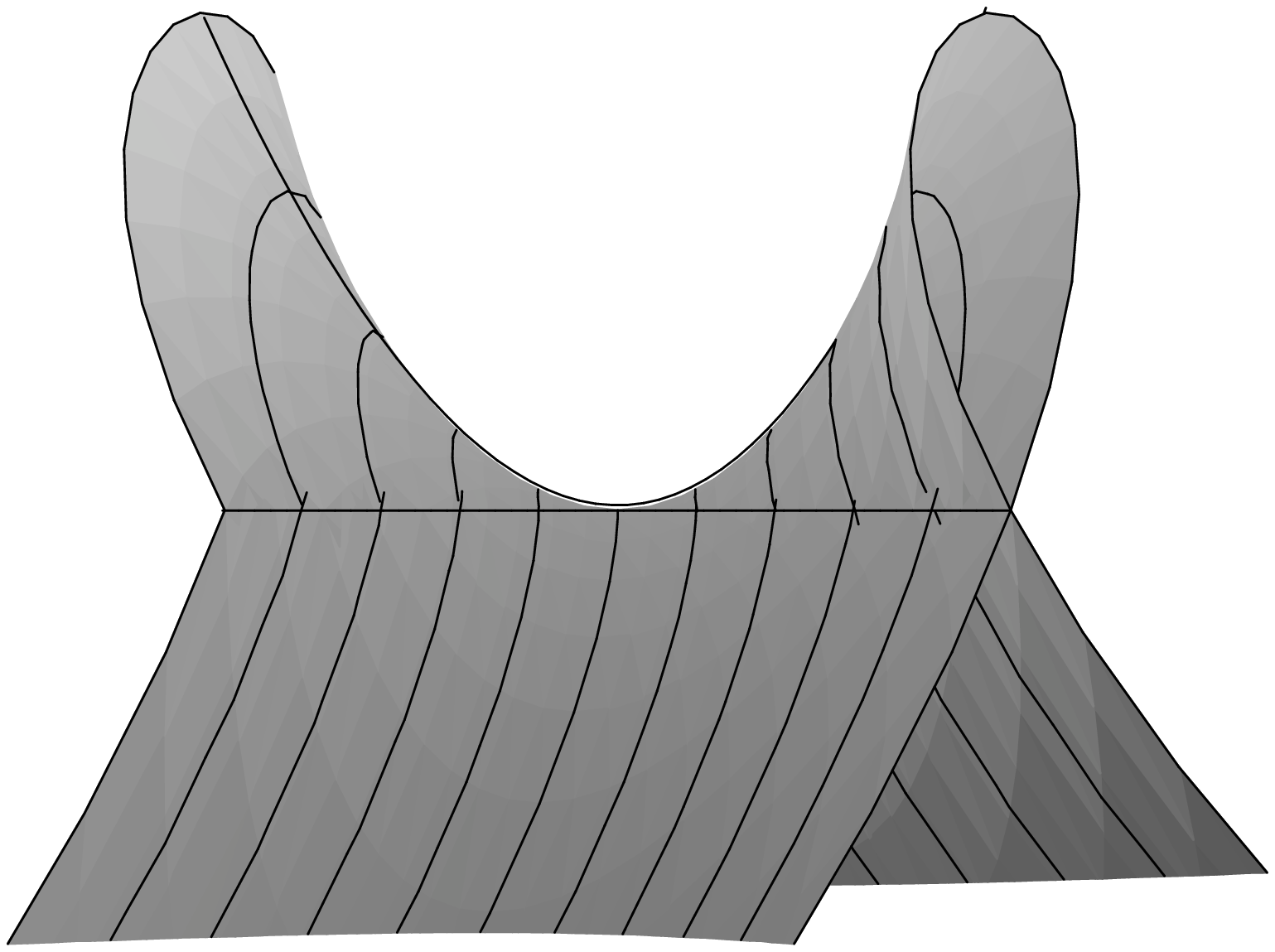,angle=90,width=9cm}}
 \put(6,0){\epsfig{file=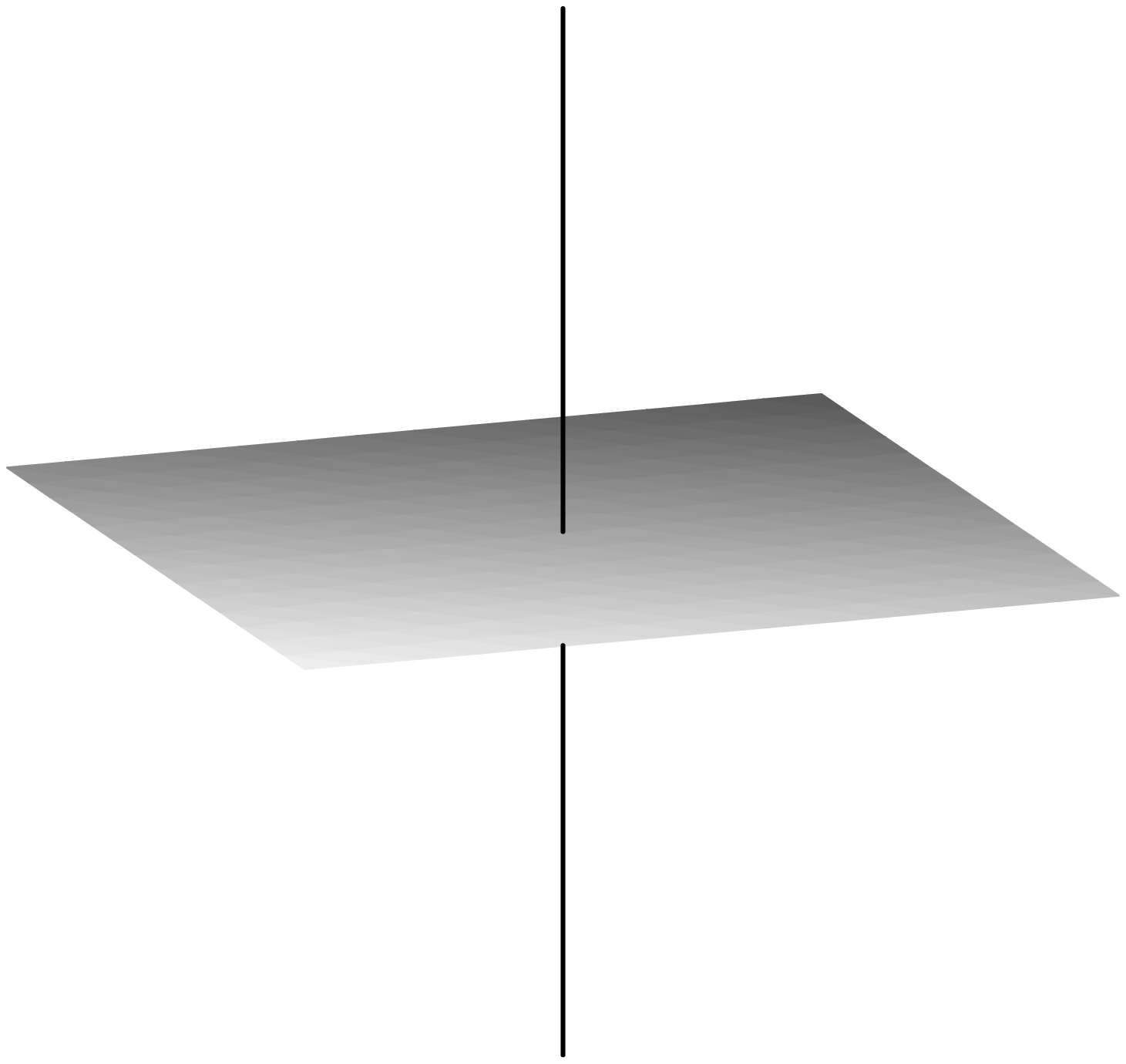,height=7.5cm,width=7.5cm}}
 \put(1.6,0.8){Example 1) the hypersurface  }
 \put(2.3,0.3){$V(x^2 + y^3- t^2 y^2)$}
 \put(8.6,0.8){Example 2) the variety}
 \put(9.5,0.3){$V(xz, yz)$}
\end{picture}

\unitlength1cm
\newsavebox{\ninepoints}
\savebox{\ninepoints}(4.5,4.5)[bl]{
  \put(0,0){
    \begin{texdraw}
      \drawdim{cm}
      \setunitscale 1.4
      \setgray 0
      \move (0 0) \fcir f:0 r:0.1
      \move (-1 -1) \fcir f:0 r:0.04
      \move (-1 0) \fcir f:0 r:0.04
      \move (-1 1) \fcir f:0 r:0.04
      \move (0 -1) \fcir f:0 r:0.04
      \move (0 1) \fcir f:0 r:0.04
      \move (1 -1) \fcir f:0 r:0.04
      \move (1 0) \fcir f:0 r:0.04
      \move (1 1) \fcir f:0 r:0.04
    \end{texdraw}
    }
  }

\unitlength1cm
\begin{picture}(13,6)
   \put(0,0){\epsfig{file=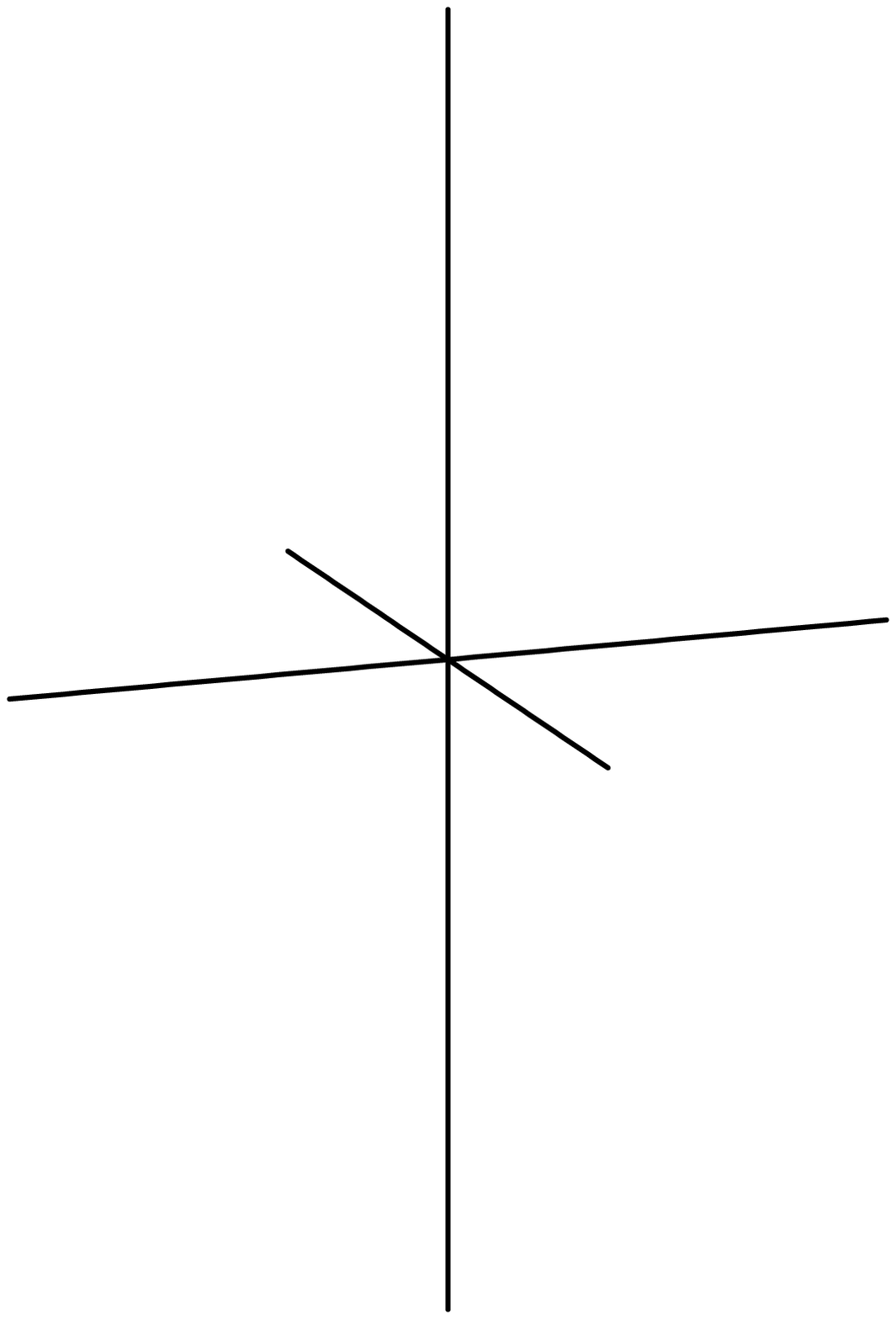,height=6.75cm,width=7.5cm}}
   \put(8,2){\usebox{\ninepoints}}
   \put(1.4,0.8){Example 3) the space curve}
   \put(2.3,0.3){$V(xy, xz, yz)$}
   \put(7.8,0.8){Example 4) the set of points}
   \put(6.7,0.3){$V(y^4-y^2,xy^3-xy,x^3y-xy,x^4-x^2)$}
\end{picture}
\medskip

\item \textit{A more geometric question is the following.   
     Let $V(I') \subset V(I)$ be a subvariety.   How can we describe $V(I)
     \smallsetminus V(I')$?  Algebraically, this amounts to finding generators for the
    ideal quotient
\[
I : I' = \{f \in K[x] \mid f I' \subset I\}.
\]
}
(The same definition applies if $I, I^\prime$ are submodules of $K[x]^k$.)

Geometrically, $V(I : I')$ is the smallest variety containing $V(I)
\smallsetminus V(I')$ which is the (Zariski) closure of $V(I) \smallsetminus
V(I')$.

In Example 2) we have $\langle xz, yz\rangle : \langle x,y\rangle = z$ and in
Example 3) $\langle xy, xz, yz\rangle : \langle x,y\rangle = \langle
z,xy\rangle$, which gives, in both cases, equations for the complement of the
$z$--axis $x = y = 0$.   In Example 4) we get $I : \langle x,y\rangle^2 =
\langle y(y^2-1), x(x^2-1), (x^2-1) (y^2-1)\rangle$ which is the zero set of
the eight points $V(I)$ with centre removed.

\item \textit{Geometrically important is the projection of a variety $V(I)
    \subset K^n$ into a linear subspace $K^{n-r}$.   Given generators $f_1,
    \dots, f_k$ of $I$, we want to find generators for the (closure of the)
    image of $V(I)$ in $K^{n-r} = \{x|x_1 = \dots = x_r = 0\}$.   The image is
    defined by the ideal $I \cap K[x_{r+1}, \dots, x_n]$ and finding
    generators for this intersection is known as eliminating $x_1, \dots, x_r$
    from $f_1, \dots, f_k$.}

Projecting the varieties of Examples 1) -- 3) to the $(x,y)$--plane is, in the first
two cases, surjective and in the third case it gives the two coordinate axes
in the $(x,y)$--plane.   This corresponds to the fact that the intersection
with $K[x,y]$ of the first two ideals is 0, while the third one is $xy$.

Projecting the 9 points of Example 4) to the $x$--axis we get, by eliminating
$y$, the polynomial $x^2 (x-1) (x+1)$, describing the three image points.
From a set theoretical point of view this is nice, however it is not
satisfactory if we wish to count multiplicities.  For example, the two border
points are the image of three points each, hence should appear with
multiplicity three.  That this is not the case can be explained by the fact
that elimination computes the annihilator ideal of $K[x,y]/I$ considered as
$K[x]$--module (and not the Fitting ideal).  This is related to the
well--known fact that elimination is not compatible with base change.

\item Another problem is related to the Riemann singularity removable theorem,
  which states that a function on a complex manifold, which is holomorphic and
  bounded outside a sub-variety of codimension 1, is actually holomorphic
  everywhere.  This is well--known for open subsets of $\C$, but in higher
  dimension there exists a second singularity removable theorem, which states
  that a function, which is holomorphic outside a sub-variety of codimension 2
  (no assumption on boundedness), is holomorphic everywhere.

For singular complex varieties this is not true in general, but those for which
the two removable theorems hold are called \textbf{normal}.   Moreover, each
reduced variety has a normalisation and there is a morphism with finite fibres
from the normalisation to the variety, which is an isomorphism outside the
singular locus.

\textit{The problem is, given a variety $V(I) \subset K^n$, find a normal
  variety $V(J) \subset K^m$ and a polynomial map $K^m \lra K^n$ inducing the
  normalisation map $V(J) \lra V(I)$.}

\textit{The problem can be reduced to irreducible varieties (but need not be,
  as we shall see) and then the equivalent algebraic problem is to find the
  normalisation of $K[x_1, \dots, x_n]/I$, that is the integral closure of
  $K[x]/I$ in the quotient field of $K[x]/I$ and present this ring as an
  affine ring $K[y_1, \dots, y_m]/J$ for some $m$ and $J$.}

For Examples 1) -- 4) it can be shown that the
normalisation of the first three varieties is smooth, the last two are the
disjoint union of the (smooth) components.  The corresponding rings are
$K[x_1, x_2],\; K[x_1, x_2] \oplus K[x_3],\; K[x_1] \oplus K[x_2] \oplus
K[x_3]$.  The fourth example has no normalisation as it is not reduced.

\textit{A related problem is to find, for a non--normal variety $V$, an ideal
  $H$ such that $V(H)$ is the non--normal locus of $V$.   The normalisation
  algorithm described below solves also this problem.}

In the examples, the non--normal locus is equal to the singular locus.

\item The significance of \textbf{singularities} appears not only in the
  normalisation problem.   The study of singularities is also called
  \textbf{local algebraic geometry} and belongs to the basic tasks of
  algebraic geometry.   Nowadays, singularity theory is a whole subject on its
  own.

A singularity of a variety is a point which has no neighbourhood in which the
Jacobian matrix of the generators has constant rank.

In Example 1) the whole $t$--axis is singular, in the three other
examples only the origin.

\textit{One task is to compute generators for the ideal of the singular locus,
  which is itself a variety.  This is just done by computing sub-determinants
  of the Jacobian matrix, if there are no components of different dimensions.
  In general, however, we need, additionally, to compute either the
  equidimensional part and ideal quotients or a primary decomposition.}

In Examples 1) -- 4), the singular locus is given by $\langle x,y\rangle,\;
\langle x,y,z\rangle,\;\langle x,y,z\rangle,\; \langle x,y\rangle^2$,
respectively.

\item \textit{Studying a variety $V(I)$, $I = (f_1, \dots, f_k)$, locally at a
    singular point, say the origin of $K^n$, means studying the ideal
    $IK[x]_{\langle x\rangle}$ generated by $I$ in the local ring
\[
K[x]_{\langle x\rangle} = \left\{\dfrac{f}{g} \mid f,g \in K[x], g \not\in \langle x_1,
  \dots, x_n\rangle\right\}.
\]
In this local ring the polynomials $g$ with $g(0) \not= 0$ are units and
$K[x]$ is a subring of $K[x]_{\langle x\rangle}$.}

\textit{Now all the problems we considered above can be formulated for ideals
  in $K[x]_{\langle x\rangle}$ and modules over $K[x]_{\langle x \rangle}$
  instead of $K[x]$.}

The geometric problems should be interpreted as properties of the
variety in a neighbourhood of the origin, or more generally, the given point.
\end{enumerate}

It should not be surprising that all the above problems have algorithmic and
computational solutions, which use, at some place, Gr\"obner basis methods.
Moreover, algorithms for most of these have been implemented quite efficiently
in several computer algebra systems, such as CoCoA, cf.\ \ccite{CNR},
Macaulay2, cf.\ \ccite{GS} and SINGULAR, cf.\ \ccite{GPS}, the latter also
being able to handle, in addition, local questions systematically.  

The most complicated problem is the primary decomposition, the latest
achievement is the normalisation, both being implemented in SINGULAR.

At first glance, it seems that computation in the localisation $K[x]_{\langle
  x\rangle}$ requires computation with rational functions.  It is an important
fact that this is not necessary, but that basically the same algorithms which
were developed for $K[x]$ can be used for $K[x]_{\langle x\rangle}$.  This is
achieved by the choice of a special ordering on the monomials of $K[x]$ where,
loosely speaking, the monomials of lower degree are considered to be bigger.

However, such orderings are no longer well--orderings and the classical
Buchberger algorithm would not terminate.  Mora discovered, cf.\ \ccite{Mo},
that a different normal form algorithm, or, equivalently, a different division
with remainders, leads to termination.  Thus, Buchberger's algorithm with
Mora's normal form is able to compute in $K[x]_{\langle x\rangle}$ without
denominators.

Several algorithms for $K[x]$ use elimination of (some auxiliary extra)
variables.  But variables to be eliminated have, necessarily, to be
well--ordered.  Hence, to be able to apply the full power of Gr\"obner basis
methods also for the local ring $K[x]_{\langle x\rangle}$, we need mixed
orders, where the monomial ordering restricted to some variables is not a
well--ordering, while restricted to other variables it is.  In \ccite{GP1} and
\ccite{Getal}, the authors described a modification of Mora's normal form,
which terminates for mixed ordering, and more generally, for any monomial
ordering which is compatible with the natural semigroup structure.

\section{Some global algorithms}

Having mentioned some geometric problems, I shall now illustrate two
algorithms related to these problems:  primary decomposition and
normalisation. 

\subsection{primary decomposition}

Any ideal $I \subset R$ in a Noetherian ring can be written as $I =
\overset{r}{\underset{i=1}{\cap}} q_i$ with $q_i$ primary ideals (that is,
$q_i \not= R$ and $gf \in q_i$ implies $g \in q_i$ or $f^p \in q_i$ for some
$p > 0$).

This generalises the unique factorisation (valid in factorial rings) $f =
f_1^{p_1} \cdot \ldots \cdot f_r^{p_r}$ with $f_i$ irreducible, from elements to
ideals.  In $K[x]$ we have both, unique factorisation and primary
decomposition and any algorithm for primary decomposition needs factorisation
(because a primary decomposition of a principal ideal $I = \langle f \rangle$
is equivalent to a factorisation of $f$).

In contrast to factorisation, primary decomposition is, in general, not unique,
even if we consider minimal decompositions, that is, the associated primes
$p_i = \sqrt{q_i}$ are all distinct and none of the $q_i$ can be omitted in
the intersection.   However, the minimal (or isolated) primes, that is, the
minimal elements of $\Ass(I) = \{p_1, \dots, p_r\}$ with regard to inclusion, are
uniquely determined.   The minimal primes are the only ``geometrically
visible'' primes in the sense that
\[
V(I) = \bigcup_{p_j \in \minAss(I)} V(p_j)
\]
is the decomposition of $V(I)$ into irreducible components.  A non--minimal
associated prime $p_i \not\in \minAss(I)$ is called embedded, because there exists a $p_j
\in \minAss(I),\; p_j \subset p_i$. This means geometrically $V(p_i) \subset
V(p_j)$, that is, the irreducible component of $V(I)$ corresponding to $p_i$
is embedded in some bigger irreducible component.

As an example we compute the primary decomposition of the ideal $I = \langle
x^2y^3-x^3yz,\; y^2z - xz^2\rangle$ in SINGULAR, the output
being slightly changed in order to save space.
\begin{verbatim}
LIB "primdec.lib";              //calling library for primary decomposition
ring R  = 0,(x,y,z),dp;
ideal I = x2y3-x3yz,y2z-xz2;
primdecGTZ(I);
==> [1]: [1]:               [2]: [1]:               [3]: [1]:
           _[1]=-y2+xz              _[1]=z2                 _[1]=z
         [2]:                       _[2]=y                  _[2]=x2
           _[1]=-y2+xz           [2]:                    [2]:
                                    _[1]=z                  _[1]=z
                                    _[2]=y                  _[2]=x
\end{verbatim}
The result is a list of three pairs of ideals (for each pair, the first ideal
is the primary component, the second ideal the corresponding prime component).
The second prime component [2] : [2] is embedded in the first [1] : [2].  The
first primary component [1] : [1] is already prime, the other two are not.

Hence, $I = (y^2 - xz) \cap (y,z^2) \cap (x^2,z)$ and we obtain:
\[
V(I) = \{y^2 - xz = 0\} \cup \underset{\text{\small (embedded
    component)}}{\{y=z^2=0\}} \cup \{x^2 = z = 0\}
\]
\begin{center}
\unitlength1cm
\begin{picture}(13,5)
\put(0,0){\includegraphics[clip,height=4cm,width=4cm]{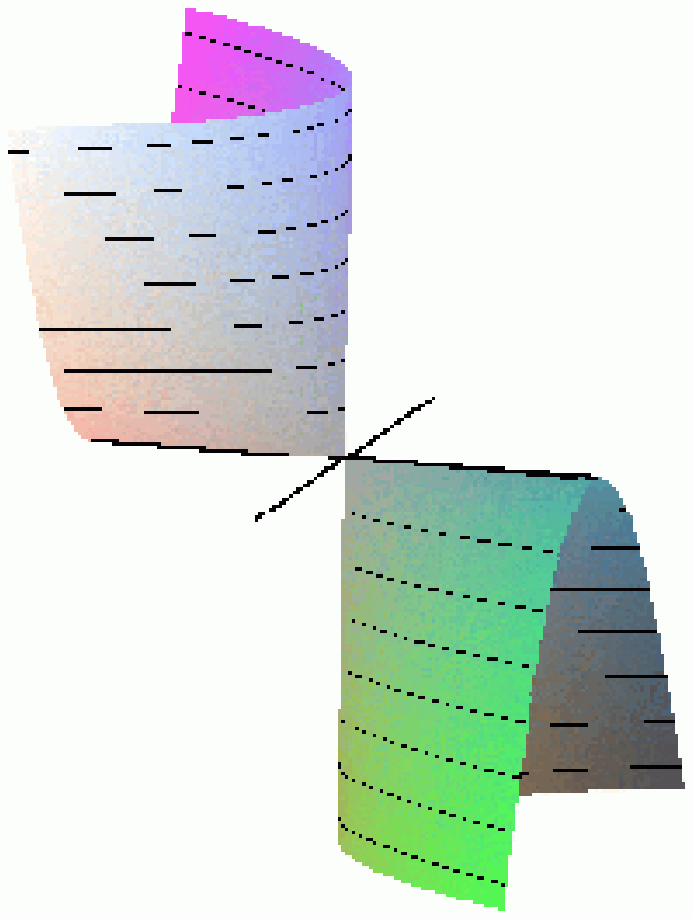}}
\put(3.5,2.5){$=$}
\put(3.8,0){\includegraphics[clip,height=4cm,width=4cm]{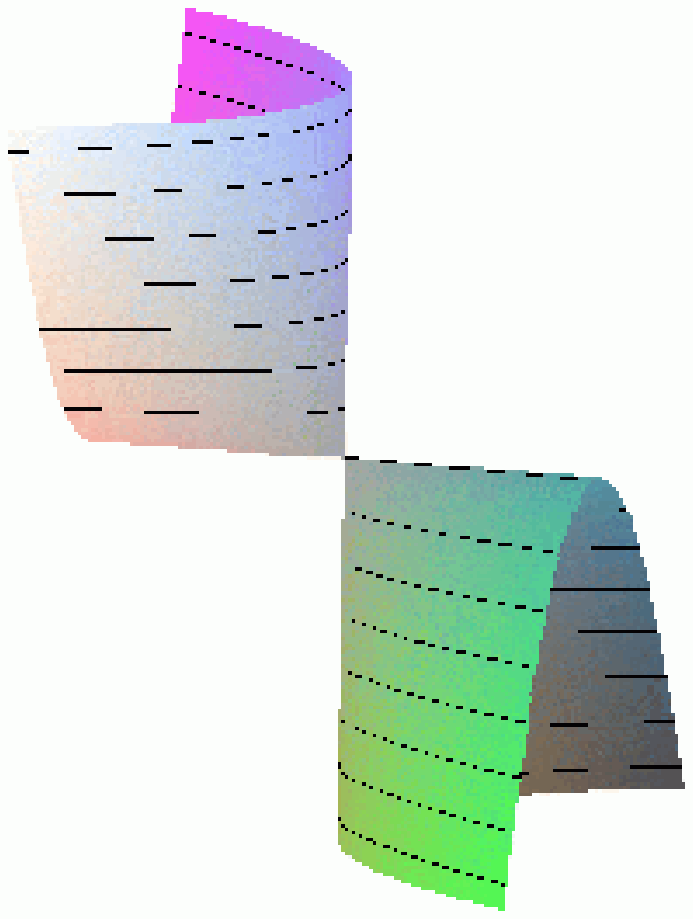}}
\put(7.5,2.5){$\cup$}
\put(8,1.5){\includegraphics[clip,height=2cm,width=2cm]{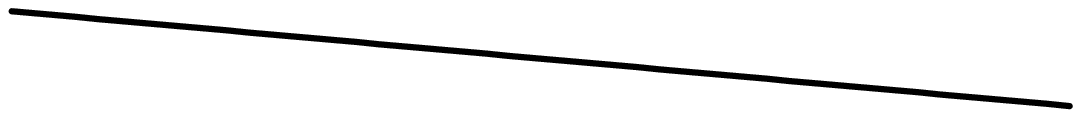}}
\put(11,2.5){$\cup$}
\put(11,1.5){\includegraphics[clip,height=2cm,width=2cm]{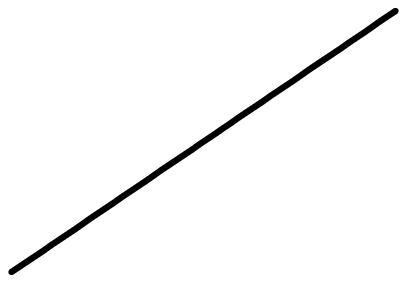}}
\end{picture}

Primary decomposition
\end{center}

All known algorithms for primary decompositions in $K[x]$ are quite involved
and use many different sub--algorithms from various parts of computer
algebra, in particular Gr\"obner bases, resp.\ characteristic sets, and
multivariate polynomial factorisation over some (algebraic or transcendental)
extension of the field $K$.  For an efficient implementation which can treat
examples of interest in algebraic geometry, a lot of extra small additional
algorithms have to be used.  In particular one should use ``easy'' splitting as
soon and as often as possible, see \ccite{DGP}.

In SINGULAR the algorithms of \ccite{GTZ} (which was the first practical and
general primary decomposition algorithm), the recent algorithm of \ccite{SY}
and some of the homological algebra algorithms for primary decomposition of
\ccite{EHV} have been implemented.  For detailed and improved versions of
these algorithms, together with extensive comparisons, see \ccite{DGP}.
\medskip

Here are some major ingredients for primary decomposition:

\begin{enumerate}
\item Reduction to zero--dimensional primary decomposition (GTZ);\\[0.5ex]
maximal independent sets;\\
ideal quotient, saturation, intersection.
\newpage
 
\item Zero--dimensional primary decomposition (GTZ);\\[0.5ex]
lexicographical Gr\"obner basis;\\
factorisation of multivariate polynomials;\\
generic change of variables;\\
primitive element computation.
\end{enumerate}

Related algorithms:
\begin{enumerate}
\item Computation of the radical;\\[0.5ex]
square--free part of univariate polynomials;\\
find (random) regular sequences (EHV).

\item Computation of the equidimensional part (EHV);\\[0.5ex]
Ext--annihilators;\\
ideal quotients, saturation and intersection.
\end{enumerate}

To see how homological algebra comes into play, let us compute the
equidimensional part of $V(I)$, that is, the union of all maximal dimensional
components of $V(I)$, or, algebraically, the intersection of all minimal
primes.  Following \ccite{EHV}, we can calculate the equidimensional part of a
variety via Ext--groups:

If $c = \codim_{K[x]}(I)$, then the equidimensional part of $I$ is the
annihilator ideal of the module Ext$^c_{K[x]} (K[x]/I, K[x])$ by \ccite{EHV}.

For example, the equidimensional part of $V = \{xz = yz = 0\}$ is given by the
ideal $\langle z \rangle = \ann\bigl(\Ext^1 (K[x,y,z]/\langle xz,yz\rangle,
K[x,y,z])\bigr)$.   
\medskip

Using SINGULAR, we obtain this via:

\begin{verbatim}
LIB "homolog.lib";  
ring  r  = 0,(x,y,z),dp;
ideal I  = xz, yz;
module M = Ext_R(1,I);
quotient(M,freemodule(nrows(M)));
==> _[1] =z
\end{verbatim}
\vspace{-4cm}

\begin{center}
  \unitlength1cm
  \begin{picture}(13,6)
\put(7,0){\epsfig{file=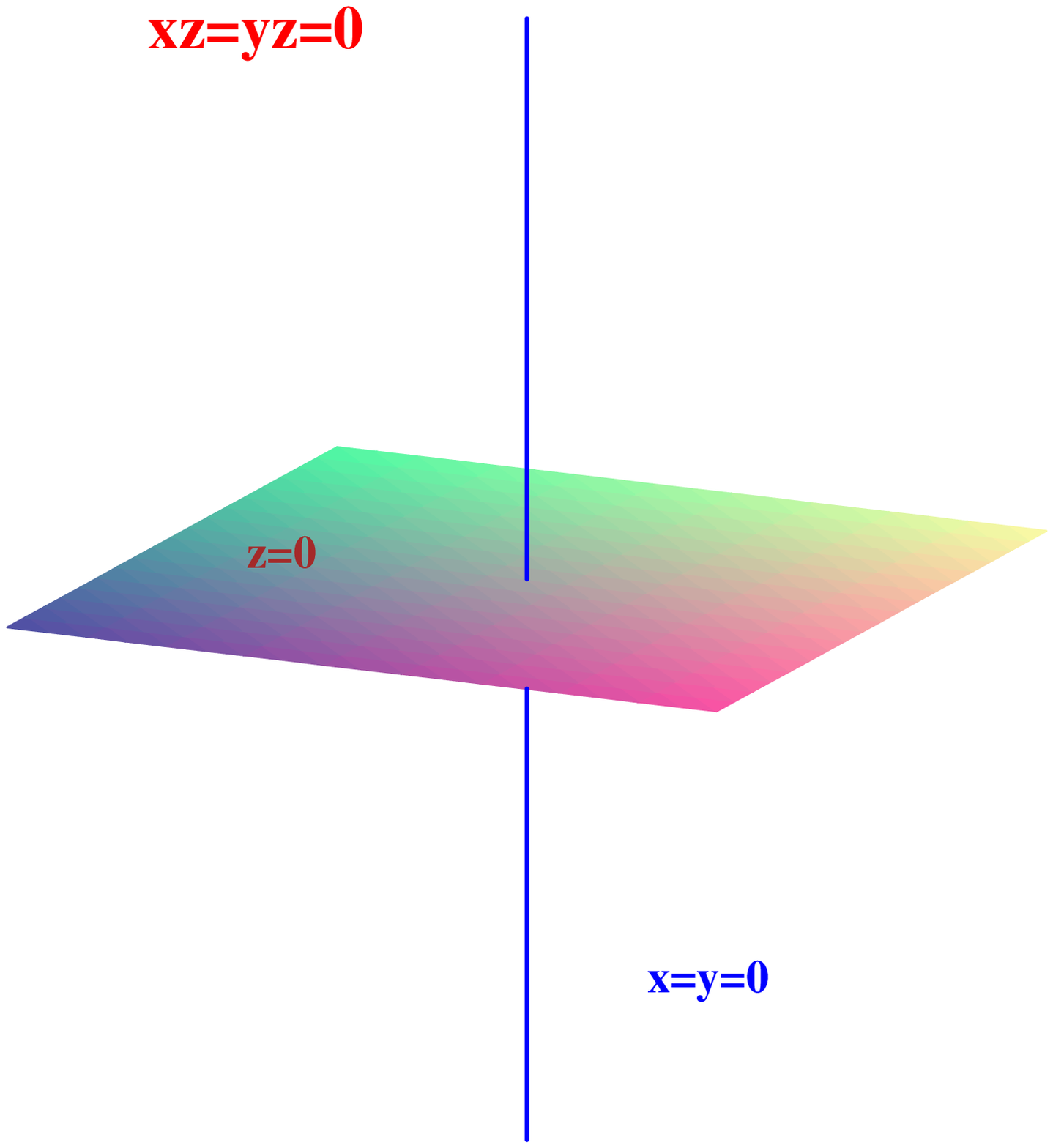,height=7cm,width=7cm}}
  \end{picture}
\end{center}
\vspace{-1.8cm}

Note that \verb?module M = Ext_R(i,I)? computes a presentation matrix of
$\Ext^i(R/I,R)$.  Hence, identifying a matrix with its column space in the
free module of rank equal to the number of rows, $\Ext^1(R/I,R) = R^n/M$ with
$R^n =$ \verb?freemodule(nrows(M))? and, therefore,
$\Ann\bigl(\Ext^1(R/I,R)\bigr) = M : R^n =$
\verb?quotient(M,freemodule(nrows(M)))?.

Above, we used the procedure \verb?Ext_R(-,-)? from \verb?homolog.lib?.  Below
we show that the Ext groups can easily be computed directly in a system which
offers free resolutions, respectively syzygies, transposition of matrices and
presentations of sub-quotients of a free module (\verb?modulo? in SINGULAR).
Indeed, the Ext--annihilator can be computed more directly (and faster)
without computing the Ext group itself:

\noindent Take a free resolution of $R/I$:
\[
0 \longleftarrow R/I \longleftarrow R \longleftarrow R^{n_1}
\longleftarrow \cdots.
\]
Then consider the dual sequence:
\[
0 \lra \Hom(R,R) \overset{d^0}{\lra} \Hom(R^{n_1}, R) \overset{d^1}{\lra}
\cdots .
\]
This leads to:

$\Ext^i(R/I,R) = \Ker(d^i)/\im(d^{i-1}) \text{ and }
\Ann\bigl(\Ext^i(R/I,R)\bigr) = \im(d^{i-1}) : \Ker(d^i)$.
\newpage

The corresponding SINGULAR commands are:
\vspace{-0.3cm}

\begin{verbatim}
int i = 1;
resolution L = res(I,i+1);
module Im    = transpose(L[i]);
module Ker   = syz(transpose(L[i+1]));
module ext   = modulo(Ker,Im);           //the Ext group
ideal ann    = quotient(Im,Ker);         //the Ext-annihilator
\end{verbatim}

Since the resolution can be computed by iterated syzygy computation, this is a
beautiful example of geometric use of syzygies.  However, the algorithm is not
at all obvious, but based on the non--trivial theorem of Eisenbud, Huneke and
Vasconcelos.

\subsection{Normalisation}

Another important algorithm is the normalisation of $K[x]/I$ where $I$ is a
radical ideal.  It can be used as a step in the primary decomposition, as
proposed in \ccite{EHV}, but is also of independent interest.  Several
algorithms have been proposed, especially by \ccite{Se}, \ccite{Sto},
\ccite{GT}, \ccite{Va1}.  It had escaped the computer algebra community,
however, that \ccite{GR} had given a constructive proof for the ideal of the
non--normal locus of a complex space.  Within this proof they provide a
normality criterion which is essentially an algorithm for computing the
normalisation, cf.\ \ccite{Jo}.  Again, to make the algorithm efficient needed
some extra work which is described in \ccite{DGJP}.  The Grauert--Remmert
algorithm is implemented in SINGULAR and seems to be the only full
implementation of the normalisation.  \smallskip

\noindent \textbf{Criterion} \cite{GR}: 
\textit{Let $R = K[x]/I$ with $I$ a radical ideal.  Let $J$ be a radical ideal
containing a non--zero divisor of $R$ such that $V(J)$ contains the
non--normal locus of $V(I)$.  Then $R$ is normal if and only if $R =
\Hom_R(J,J)$.}

For $J$ we may take any ideal so that $V(J)$ contains the singularities of
$V(I)$.   Since normalisation commutes with localisation, we obtain
\smallskip

\noindent \textbf{Corollary}: \textit{Ann$(\Hom_R(J,J)/R)$ is an ideal describing the non--normal locus of $V(I)$.}
\medskip

Now $\Hom_R(J,J)$ is a ring containing $R$ and if $R \subsetneqq
\Hom_R(J,J)=R_1$ we can continue with $R_1$ instead of $R$ and obtain an
increasing sequence of rings $R \subset R_1 \subset R_2 \subset \ldots$.

After finitely many steps the sequence becomes stationary (because the
normalisation of $R = K[x]/I$ is finite over $R$) and
we reach the normalisation of $R$ by the criterion of Grauert and Remmert.

Ingredients for the normalisation (which is a highly recursive algorithm):  
\begin{enumerate}
\item Computation of the ideal $J$ of the singular locus of the ideal $I$;
\item computation of a non--zero divisor for $J$;
\item ring structure on $\Hom(J,J)$;
\item syzygies, normal forms, ideal quotient.
\end{enumerate}

\noindent SINGULAR commands for computation of the normalisation:
\vspace{0.5cm}

\begin{minipage}[b]{4.5cm}
\begin{verbatim}
LIB "normal.lib";
ring S   = 0,(x,y,z),dp;
ideal I  = y2-x2z;
list nor = normal(I);
def R    = nor[1];
setring R;
normap;
==> normap[1]=T(1)
==> normap[2]=T(1)*T(2)
==> normap[3]=T(2)^2
\end{verbatim}
\end{minipage}
\hfill \parbox[b]{8.5cm}{
 \unitlength0.75cm
  \begin{picture}(8.5,4)
    \put(-1.5,-1){\includegraphics[clip,width=4cm]{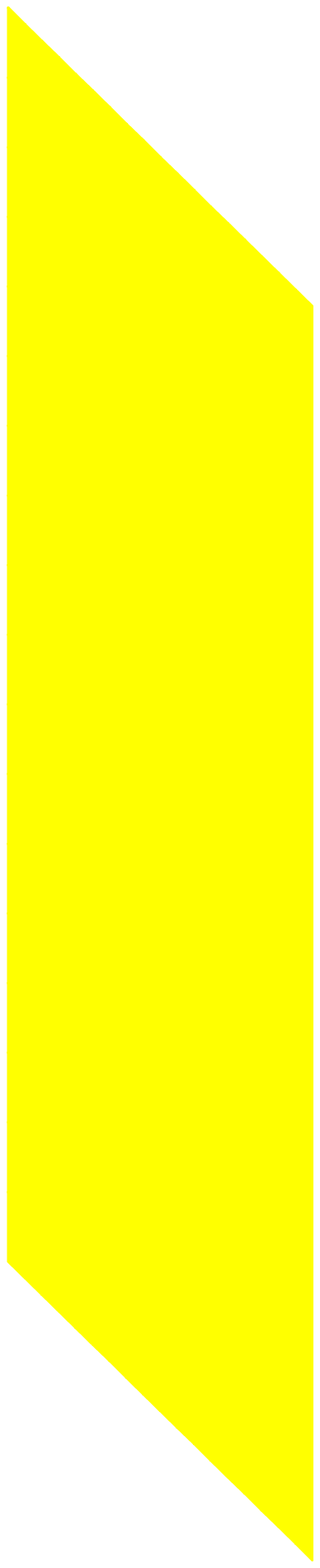}}
    \put(1.5,-0.6){\includegraphics[clip,width=3.8cm]{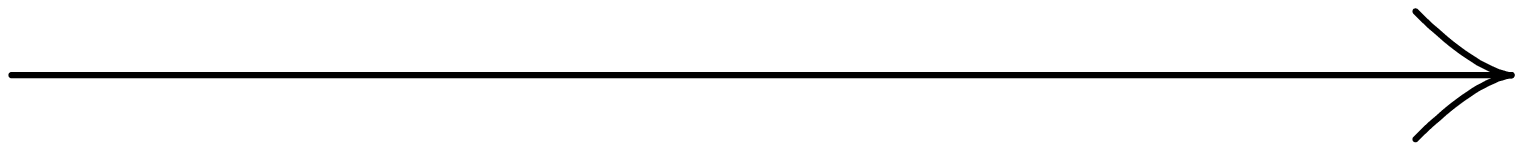}}
    \put(2,1.7){\parbox{3cm}{\begin{center}$(s,t)\mapsto(s,st,t^2)$\end{center}}}
    \put(6.5,-0.3){\includegraphics[clip,width=5cm]{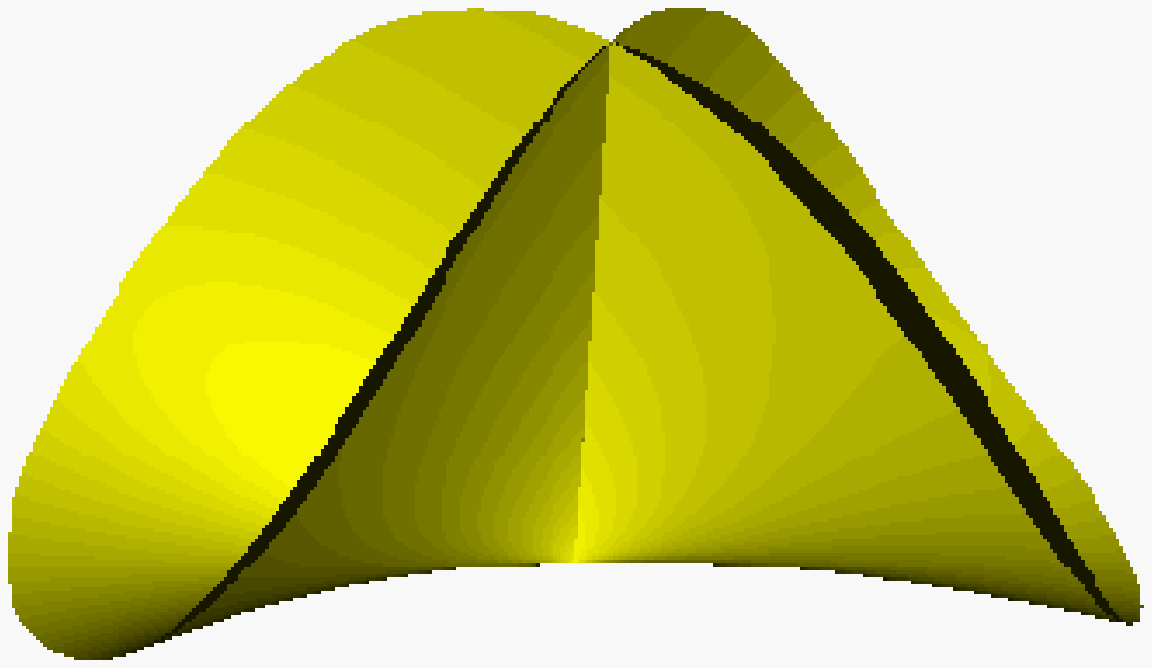}}    
  \end{picture}}
\vspace{0.3cm}

\noindent In the preceding picture, $R$, the normalisation of $S$, is
  just the polynomial ring in two variables $T(1)$ and $T(2)$.  (The
  ``handle'' of Whitney's umbrella is invisible in the parametric picture
  since it requires an imaginary parameter $t$.)
\smallskip

In several cases the normalisation of a variety is smooth (for example, the
normalisation of the discriminant of a versal deformation of an isolated
hypersurface singularity) sometimes even an affine space.   In this case, the
normalisation map provides a parametrisation of the variety.   This is the
case for the Whitney umbrella:  $V = \{y^2 - zx^2 = 0\}$.

\section{Singularities and standard bases}

A {\bf (complex) singularity} is, by definition, nothing but a complex
analytic germ $(V,0)$ together with its analytic local ring $R = \C\{x\}/I$,
where $\C\{x\}$ is the convergent power series ring in $x = x_1, \dots, x_n$.
For an arbitrary field $K$ let $R = K[[x]] /I$ for some ideal $I$
in the formal power series ring $K[[x]]$.  We call $(V,0) = (\Spec R, \fm)$ or
just $R$ a singularity ($\fm$ denotes the maximal ideal of the local ring $R$)
and write $K\langle x \rangle$ for the convergent and for the formal power
series ring if the statements hold for both.

If $I \subset K[x]$ is an ideal with $I \subset \langle x \rangle = \langle
x_1, \dots, x_n\rangle$ then the singularity of $V(I)$ at $0 \in K^n$ is,
using the above notation, $K \langle x\rangle /I \cdot K\langle x\rangle $.
However, we may also consider the local ring $K[x]_{\langle x \rangle}/I \cdot
K[x]_{\langle x \rangle}$ with $K[x]_{\langle x \rangle}$ the localisation of
$K[x]$ at $\langle x \rangle$, as the singularity of $V(I)$ at $0$.
Geometrically, for $K = \C$, the difference is the following: $\C\{x\}/I
\C\{x\}$ describes the variety $V(I)$ in an arbitrary small neighbourhood of
$0$ in the Euclidean topology while $\C[x]_{\langle x \rangle}/I
\C[x]_{\langle x \rangle}$ describes $V(I)$ in an arbitrary small
neighbourhood of $0$ in the (much coarser) Zariski topology.

At the moment, we can compute efficiently only in $K[x]_{\langle x \rangle}$
as we shall explain below.  In many cases of interest, we are happy
since invariants of $V(I)$ at $0$ can be computed in $K[x]_{\langle x \rangle}$
as well as in $K\langle x\rangle$.  There are, however, others (such as
factorisation), which are completely different in both rings.

\begin{center}
\unitlength0.8cm
\begin{picture}(16,13)
\put(5.9,12){\bf Isolated Singularities}
 \put(1,6.9){\includegraphics[clip,bb=0 150 612 592,width=5.04cm]{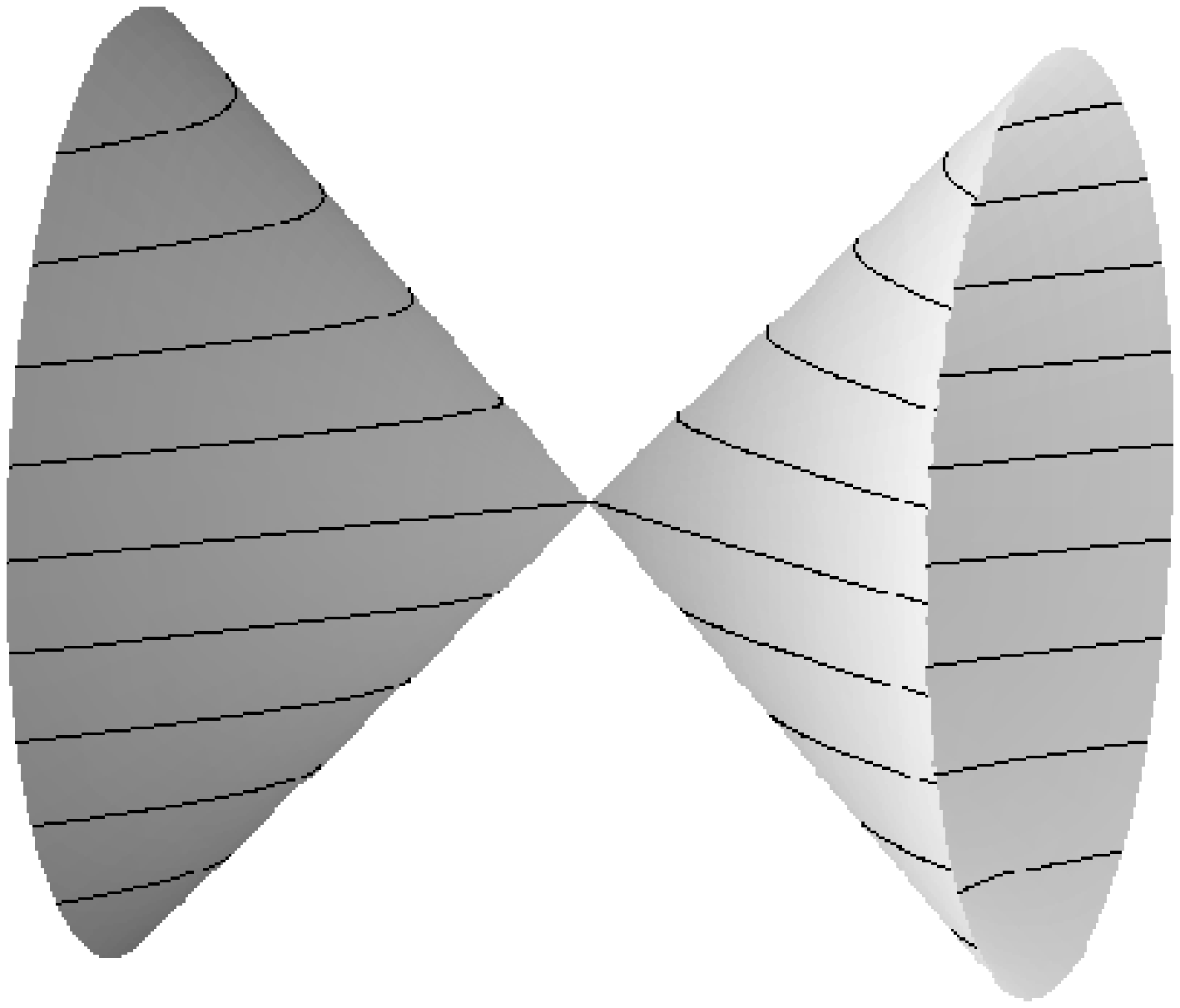}}
\put(9,6.4){\includegraphics[clip,bb=0 150 612 592,width=5.76cm]{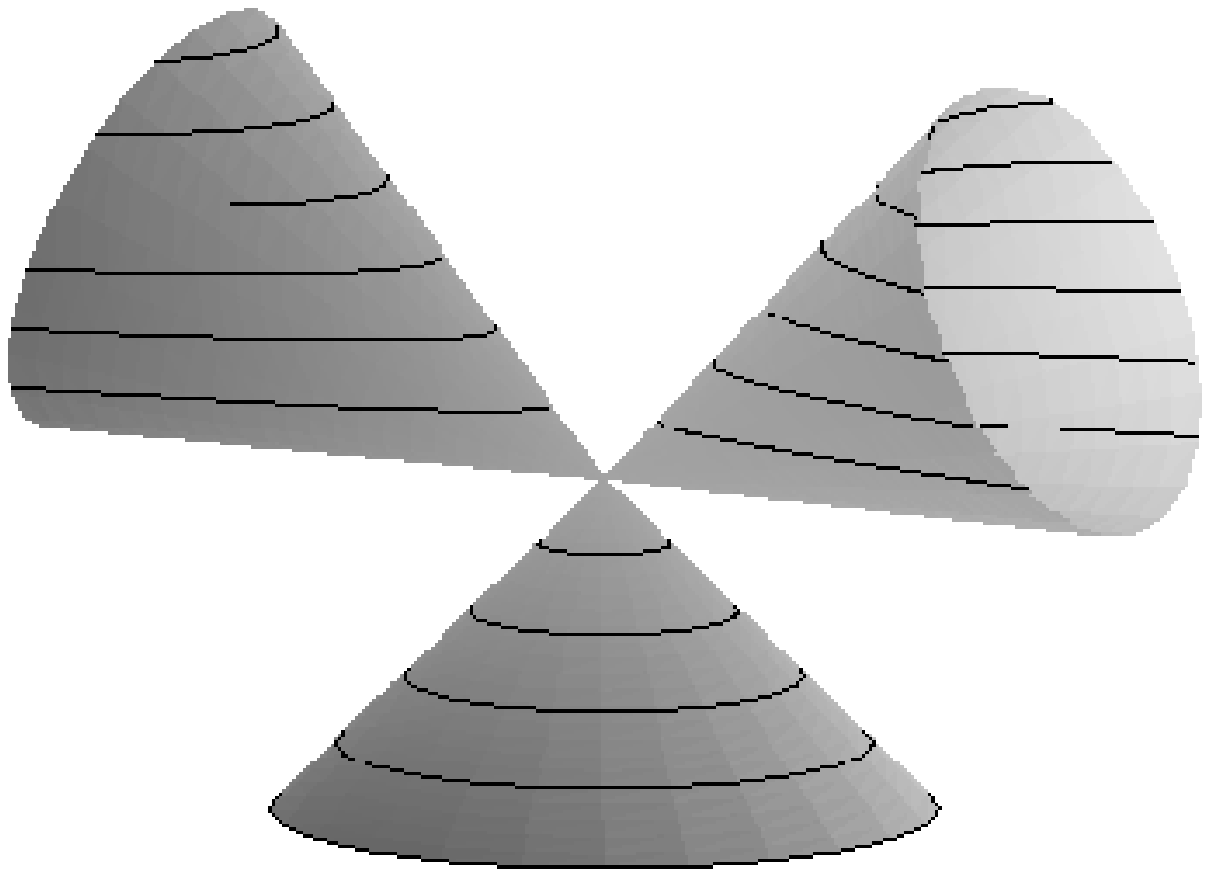}}
   \put(5.7,6){\bf Non--isolated singularities}
   \put(1,0.3){\includegraphics[clip,bb=186 261 386 531,width=3.54cm]{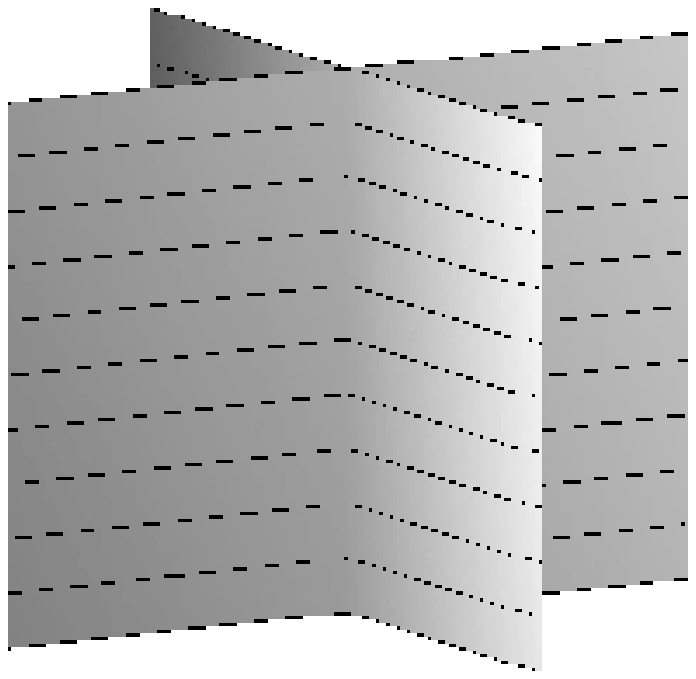}}
   \put(9.6,0.3){\includegraphics[clip,bb=0 150 612 592,width=4.32cm]{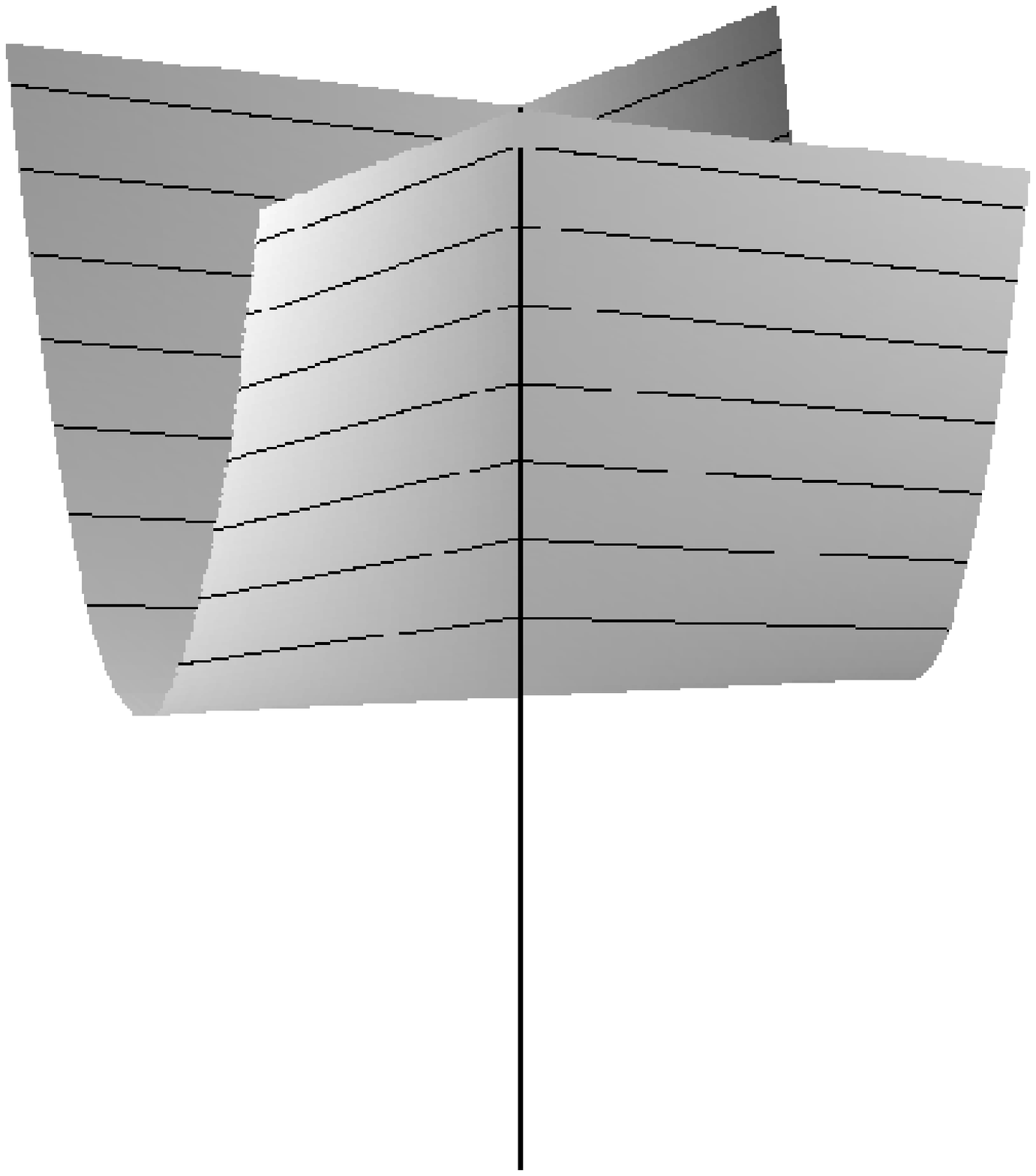}}
   \put(2.3,6.8){$A_1 : x^2 - y^2 + z^2 = 0$ }
   \put(10.5,6.8){$D_4 : z^3 - zx^2 + y^2 = 0$}
   \put(2.3,0){$A_\infty : x^2 - y^2=0$}
   \put(11,0){$D_\infty : y^2 - zx^2 = 0$}
\end{picture}
\end{center}

$(V,0)$ is called {\bf non--singular} or {\bf regular} or {\bf smooth} if
$K\langle x \rangle/I$ is isomorphic (as local ring) to a power series ring
$K\langle y_1, \dots, y_d\rangle$, or if
$K[x]_{\langle x \rangle}/I$ is a regular local ring.

By the implicit function theorem, or by the Jacobian criterion, this is
equivalent to the fact that $I$ has a system of generators $g_1, \dots,
g_{n-d}$ such that the Jacobian matrix of $g_1, \dots, g_{n-d}$ has rank $n-d$
in some neighbourhood of $0$.  $(V,0)$ is called an {\bf isolated singularity}
if there is a neighbourhood $W$ of $0$ such that $W \cap (V
\smallsetminus\{0\})$ is regular everywhere.
\medskip

In order to compute with singularities, we need the notion of standard basis
which is a generalisation of the notion of Gr\"obner basis, cf.\ 
\ccite{GP1,GP2}.

A {\bf monomial ordering} is a total order on the set of monomials
  $\{x^\alpha | \alpha \in \N^n\}$ satisfying
\[
x^\alpha > x^\beta \Rightarrow x^{\alpha + \gamma} > x^{\beta + \gamma} \text{
  for all } \alpha, \beta, \gamma \in \N^n.
\]
We call a monomial ordering $>$ \textbf{global} (resp.\ \textbf{local}, resp.\ 
\textbf{mixed}) if $x_i > 1$ for all $i$ (resp.\ $x_i< 1$ for all $i$, resp.\ 
if there exist $i,j$ so that $x_i <1$ and $x_j > 1$).  This notion is
justified by the associated ring to be defined below.  Note that $>$ is global
if and only if $>$ is a well--ordering (which is usually assumed).

Any $f \in K[x]\smallsetminus\{0\}$ can be written uniquely as $f = cx^\alpha
+ f'$, with $c \in K\smallsetminus\{0\}$ and $\alpha > \alpha'$ for any
non--zero term $c' x^{\alpha'}$ of $f'$.  We set lm$(f) = x^\alpha$, the {\bf
  leading monomial} of $f$ and $\Lc(f) = c$, the {\bf leading coefficient} of
$f$.

For a subset $G \subset K[x]$ we define the {\bf leading ideal} of $G$ as
\[
L(G) = \langle \text{ lm}(g) \;|\; g \in G\smallsetminus\{0\}\rangle_{K[x]},
\]
the ideal generated by the leading monomials in $G\smallsetminus\{0\}$.
\medskip

So far, the general case is not different to the case of a well--ordering.
 However, the following definition provides something new for non--global
 orderings:

For a  monomial ordering $>$ define the multiplicatively closed set
\[
S_>  := \{u \in K[x]\smallsetminus\{0\}\;|\;\lm (u) = 1\}
\]
and the $K$--algebra
\[
R:= \Loc K[x]  := S^{-1}_> K[x] = \{\dfrac{f}{u} \;|\; f \in K[x], u \in
S_>\},
\]
the localisation (ring of fractions) of $K[x]$ with respect to $S_>$.  We call
$\Loc K[x]$ also the {\bf ring associated to} $K[x]$ {\bf and} $>$.

Note that $K[x] \subset \Loc K[x] \subset K[x]_{\langle x\rangle}$ and $\Loc
K[x] = K[x]$ if and only if $>$ is global and $\Loc K[x] = K[x]_{\langle
  x\rangle}$ if and only if $>$ is local (which justifies the names).

Let $>$ be a fixed monomial ordering.   In order to have a short
notation, I write
\[
R := \Loc K[x] = S^{-1}_> K[x]
\]
to denote the localisation of $K[x]$ with respect to $>$.
\medskip

Let $I \subset R$ be an ideal.   A finite set $G \subset I$ is
    called a {\bf standard basis} of $I$ if and only if $L(G) = L(I)$, that
    is, for any $f \in I\smallsetminus\{0\}$ there exists a $g \in G$
    satisfying $\lm(g)|\lm(f)$.

If the ordering is a well--ordering, then a standard basis $G$ is
    called a {\bf Gr\"obner basis}.   In this case $R = K[x]$ and, hence, $G
    \subset I \subset K[x]$.
\medskip
    
Standard bases can be computed in the same way as Gr\"obner bases except that
we need a different \textit{normal form}.  This was first noticed by
\ccite{Mo} for local orderings (called tangent cone orderings by Mora) and, in
general, by \ccite{GP1,Getal}.

Let $\kg$ denote the set of all finite and ordered subsets $G \subset R$.
A map
\[
\NF : R \times \kg \to R,\; (f,G) \mapsto \NF(f|G),
\]
is called a {\bf normal form} on $R$ if, for all $f$ and $G$,

\begin{enumerate}
\item[(i)] $\NF(f|G) \not= 0 \Rightarrow \lm\bigl(\NF(f|G)\bigr) \not\in
  L(G)$,
\item[(ii)] $f - \NF(f|G) \in \langle G\rangle_R$, the ideal in $R$ generated
  by $G$.   
\end{enumerate}

 NF is called a {\bf weak normal form} if, instead of
(ii), only the following condition (ii') holds:

\begin{enumerate}
\item[(ii')] for each $f \in R$ and each $G \in \kg$ there exists a unit $u
  \in R$, so that $uf-\NF(f|G) \in \langle G\rangle_R$.
\end{enumerate}

Moreover, we need (in particular for computing syzygies) (weak) normal forms
with \textbf{standard representation}: if $G = \{g_1, \dots, g_k\}$, we can
write
\[
f - \NF(f|G) = \sum^k_{i=1} a_i g_i,\quad a_i \in R,
\]
such that $\lm\bigl(f-\NF(f|G)\bigr) \ge \lm(a_ig_i)$ for all $i$, that is, no
cancellation of bigger leading terms occurs among the $a_i g_i$.
\bigskip

Indeed, if $f$ and $G$ consist of polynomials, we can compute, in finitely
many steps, weak normal forms with standard representation such that $u$ and
$\NF(f|G)$ are polynomials and, hence, compute polynomial standard bases which
enjoy most of the properties of Gr\"obner bases.

Once we have a weak normal form with standard representation, the general
standard basis algorithm may be formalised as follows:
\smallskip

\noindent\textsc{Standardbasis}(G,NF) [\,arbitrary monomial ordering\,]

\noindent\textit{Input:}\;\; $G$ a finite and ordered set of polynomials, NF a
weak normal form with standard representation.

\noindent\textit{Output:} $S$ a finite set of polynomials which is a standard
basis of $\langle G\rangle_R$.
\smallskip

\noindent -- $S = G$;\\
-- $P = \{(f,g) \mid f,g\in S\}$;\\
-- while $(P \not= \emptyset)$\\
\hspace*{1cm} choose $(f,g) \in P$;\\
\hspace*{1cm} $P = P \smallsetminus \{(f,g)\}$;\\
\hspace*{1cm} $h = \NF(\text{spoly}(f,g)\mid S)$;\\
\hspace*{1cm} if $(h \not= 0)$\\
\hspace*{1.5cm} $P = P \cup \{(h,f) \mid f \in S\}$;\\
\hspace*{1.5cm} $S = S \cup \{h\}$;\\
-- return $S$;
\smallskip

Here spoly$(f,g) = x^{\gamma - \alpha} f - \tfrac{\lc(f)}{\lc(g)} x^{\gamma -
  \beta} g$ denotes the $\boldsymbol{s}$\textbf{--polynomial} of $f$ and $g$
where $x^\alpha = \lm(f),\; x^\beta = \lm(g),\; \gamma = \lcm(\alpha,\beta)$.
 
The algorithm terminates by Dickson's lemma or by the noetherian property of
the polynomial ring (and since NF terminates).  It is correct by Buchberger's
criterion, which generalises to non--well--orderings. 

If we use Buchberger's normal form below, in the case of a well--ordering,
\textsc{Standardbasis} ist just  Buchberger's algorithm:
\smallskip

\noindent\textsc{NFBuchberger}(f,G) [\,well--ordering\,]

\noindent\textit{Input:}\;\; $G$ a finite ordered set of polynomials, $f$ a
polynomial.

\noindent\textit{Output:} $h$ a normal form of $f$ with respect to $G$ with
standard representation.
\smallskip

\noindent -- $h = f$;\\
-- while $(h \not= 0$ and exist $g \in G$ so that $\lm(g)\mid\lm(h))$\\
\hspace*{1cm} choose any such $g$;\\
\hspace*{1cm} $h = \text{spoly}(h,g)$;\\
-- return $h$;
\smallskip

For an algorithm to compute a weak normal form in the case of an arbitrary
ordering, we refer to \ccite{GP1}.  \smallskip
\medskip

To illustrate the difference between local and global orderings, we compute
the dimension of a variety at a point and the (global) dimension of the
variety.

The \textbf{dimension} of the singularity $(V,0)$, or the dimension of $V$ at
$0$, is, by definition, the Krull dimension of the analytic local ring
$\ko_{V,0} = K\langle x \rangle/I$, which is the same as the Krull dimension
of the algebraic local ring $K[x]_{\langle x\rangle}/I$ in case $I = \langle f_1,
\dots, f_k\rangle$ is generated by polynomials, which follows easily from the
theory of dimensions by Hilbert--Samuel series.

Using this fact, we can compute $\dim(V,0)$ by computing a standard basis of
the ideal $\langle f_1, \dots, f_k\rangle$ generated in $\Loc K[x]$ with respect to any
\textit{local}\/ monomial ordering on $K[x]$.   The dimension is equal to the
dimension of the corresponding monomial ideal (which is a combinatorial
problem).

For example, the dimension of the affine variety $V = V(yx-y, zx-z)$ is 2 but
the dimension of the singularity $(V,0)$ (that is, the dimension of $V$ at the
point $0$) is 1:

  \begin{center}
      \epsfig{file=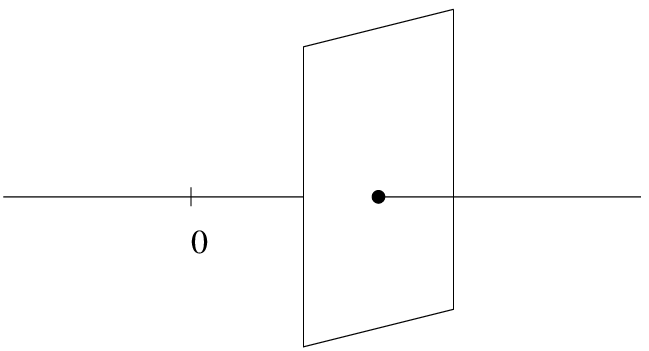}\\[2.0ex]
$V : y(x-1) = z(x-1) = 0$,\\
$ \dim(V,0) = 1,\; \dim V = 2$
  \end{center}
  
  Using SINGULAR we compute first the global dimension with the degree reverse
  lexicographical ordering denoted by dp and then the local dimension at $0$
  using the negative degree reverse lexicographical ordering denoted by ds.
  Note that in the local ring $K[x,y]_{\langle x,y\rangle}$ (represented by
  the ordering ds) $x-1$ is a unit.

\begin{verbatim}
ring R   = 0,(x,y,z),dp;  //global ring
ideal i  = yx-y,zx-z;
ideal si = groebner(i);
si;
==> si[1]=xz-z,           //leading ideal of i is <xz,xy>  
==> si[2]=xy-y
dim(si);
==> 2                     //global dimension = dim R/<xz,xy>

ring r   = 0,(x,y,z),ds;  //local ring
ideal i  = yx-y,zx-z;
ideal si = groebner(i);
si;
==> si[1]=y               //leading ideal of i is <y,z>
==> si[2]=z
dim(si);
==> 1                     //local dimension = dim r/<y,z>  
\end{verbatim}

\section{Some local algorithms}

I describe here three algorithms which use, in an essential way, standard
bases for local rings:  classification of singularities, deformations and the
monodromy.

\subsection{Classification of singularities}

In a tremendous work, V.~I.~Arnold started, in the late sixties, the
classification of hypersurface singularities up to right equivalence.  Here
$f$ and $g \in K\langle x_1, \dots, x_n\rangle$ are called right equivalent if
they coincide up to analytic coordinate transformation, that is, if there
exists a local $K$--algebra automorphism $\varphi$ of $K\langle x \rangle$
such that $f = \varphi(g)$.  His work culminated in impressive lists of
normal forms of singularities and, moreover, in a determinator for
singularities which allows the determination of the normal form for a given
power series ([AGV, II.16]).  This work of Arnold has found numerous
applications in various areas of mathematics, including singularity theory,
algebraic geometry, differential geometry, differential equations, Lie group
theory and theoretical physics.  The work of Arnold was continued by
C.T.C.~Wall and others, cf.\ \ccite{Wa,GKr}.

Most prominent is the list of ADE or simple or Kleinian singularities, which
have appeared in surprisingly different areas of mathematics, and still today,
new connections of these singularities to other areas are being discovered
(cf.\ \ccite{Gre2} for a survey).   Here is the list of ADE singularities (the
names come from their relation to the simple Lie groups of type A, D and E).
\[
\begin{array}{cclc}
 A_k & : & x_1^{k+1} + x_2^2 + x_3^2       + \dots + x_n^2, & k \ge 1 \\
 D_k & : & x_1(x_1^{k-2} + x_2^2) + x_3^2  + \dots + x_n^2, & k \ge 4 \\
 E_6 & : & x_1^4 + x_2^3 + x_3^2           + \dots + x_n^2,  & \\
 E_7 & : & x_2(x_1^3 + x_2^2) + x_3^2      + \dots + x_n^2,  & \\
 E_8 & : & x_1^5 + x_2^3 + x_3^2           + \dots + x_n^2.  &
\end{array}
\]

\begin{center}
   \unitlength1cm
   \begin{picture}(13,5)
     \put(0,0.5){\includegraphics[clip,height=4.5cm,width=4.5cm]{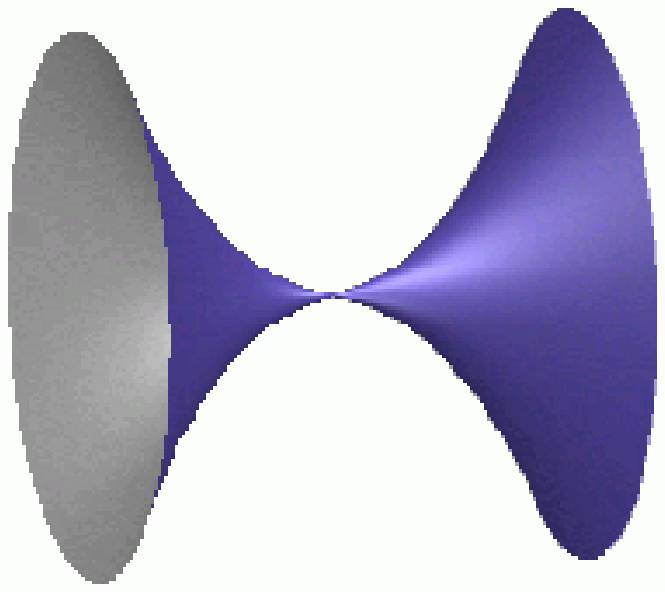}}
     \put(1,0){$A_3$-Singularity}
     \put(4.5,0.5){\includegraphics[clip,height=4.5cm,width=4.5cm]{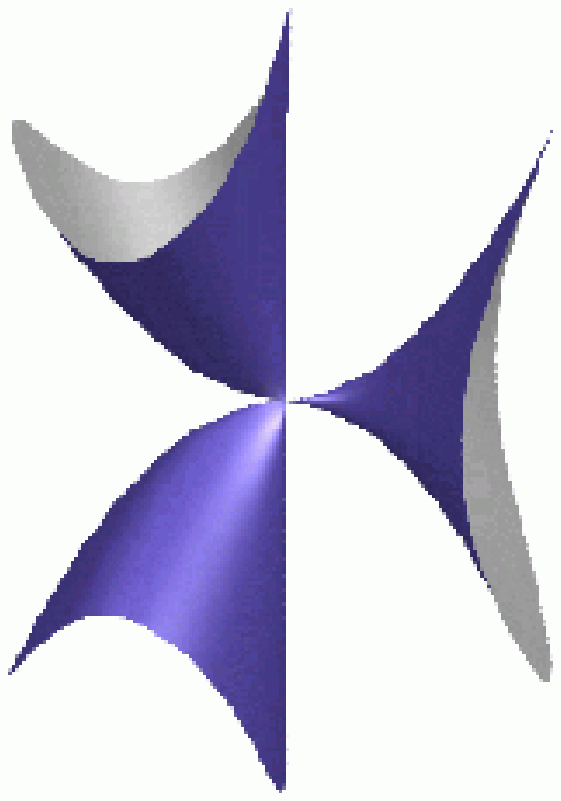}}
     \put(5.5,0){$D_6$-Singularity}
     \put(9,0.5){\includegraphics[clip,height=4.2cm,width=4.2cm]{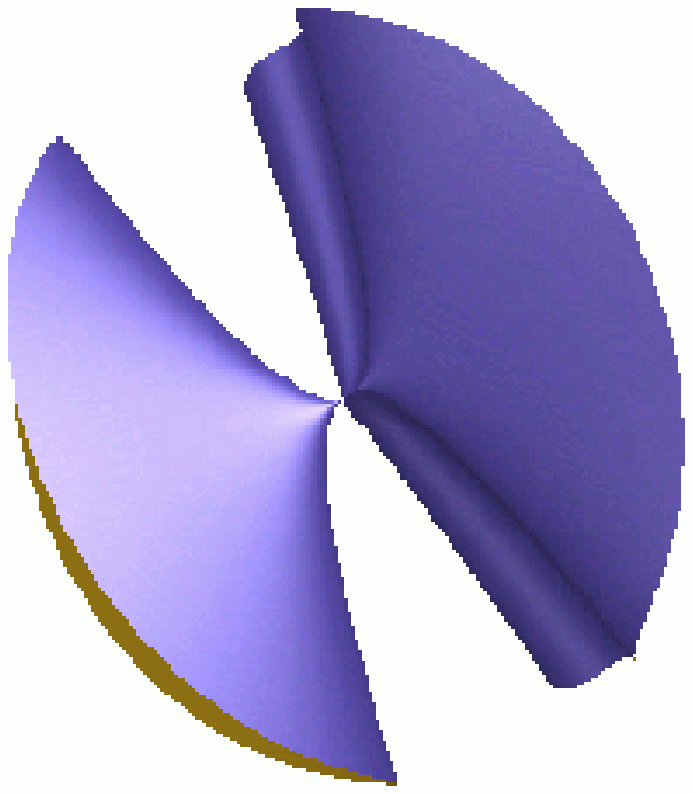}}
     \put(10,0){$E_7$-Singularity}
   \end{picture}
\end{center}

Arnold introduced the concept of ``modality'', related to Riemann's idea of
moduli, into singularity theory and classified all singularities of modality
$\le 2$ (and also of Milnor number $\le 16$).  The ADE singularities are just
the singularities of modality 0.  Singularities of modality 1 are the three
parabolic singularities:
\[
\begin{array}{lcll}
\widetilde{E}_6 = P_8 = T_{333} & : & x^3+y^3 +z^3 + axyz,& a^3 + 27 \not= 0,\\
\widetilde{E}_7 = X_g = T_{244} & : & x^4 + y^4 + ax^2y^2,& a^2 \not= 4,\\
\widetilde{E}_8 = J_{10} = T_{236} & : & x^3 + y^5 + ax^2y^2, & 4a^3 + 27
\not= 0,
\end{array}
\]
the 3--indexed series of hyperbolic singularities
\[T_{pqr} : x^p + y^q + z^r + axyz, a \not= 0, \dfrac{1}{p} + \dfrac{1}{q} +
\dfrac{1}{r} < 1
\]
and 14 exceptional families, cf.~\ccite{AGV}.

The proof of Arnold for his determinator is, to a great part constructive, and
has been partly implemented in SINGULAR, cf.\ \ccite{Kr}.  Although the whole
theory and the proofs deal with power series, everything can be reduced to
polynomial computation since we deal with isolated singularities, which are
\textbf{finitely determined}.  That is, for an isolated singularity $f$, there
exists an integer $k$ such that $f$ and $g$ are right equivalent if their
Taylor expansion coincides up to order $k$.  Therefore, knowing the
determinacy $k$ of $f$, we can replace $f$ by its Taylor polynomial up to
order $k$.

The determinacy can be estimated as the minimal $k$ such that
\[
{\fm}^{k+1} \subset {\fm}^2 \text{ jacob}(f)
\]
where ${\fm} \subset K\langle x_1, \dots, x_n\rangle$ is the maximal ideal and
jacob$(f) = \langle \partial f/\partial x_1, \dots, \partial f/\partial
x_n\rangle$.   Hence, this $k$ can be computed by computing a standard basis
of ${\fm}^2$ jacob$(f)$ and normal forms of ${\fm}^i$ with respect to this
standard basis for increasing $i$, using a local monomial
ordering. However, there is a much faster way to compute the determinacy
directly from a standard basis of ${\fm}^2 \text{ jacob}(f)$, which is basically the
``highest corner'' described in \ccite{GP1}.

An important initial step in Arnold's classification is the generalised Morse
lemma, or splitting lemma, which says that $f \circ \varphi(x_1, \dots, x_n) =
x_1^2 + \dots + x_r^2 + g(x_{r+1}, \dots, x_n)$ for some analytic coordinate
change $\varphi$ and some power series $g \in {\fm}^3$ if the rank of the
Hessian matrix of $f$ at 0 is $r$.

The determinacy allows the computation of $\varphi$ up to sufficiently high
order and a polynomial $g$ as in the theorem.   This has been
implemented in SINGULAR and is a cornerstone in classifying hypersurface
singularities. 

In the following example we use SINGULAR to get the singularity $T_{5,7,11}$
from a database A$_-$L (``Arnold's list''), make some coordinate change and
determine then the normal form of the complicated polynomial after coordinate
change.
\begin{verbatim}
LIB "classify.lib";
ring r  = 0,(x,y,z),ds;
poly f  = A_L("T[5,7,11]"); 
f;
==> xyz+x5+y7+z11
map phi = r, x+z,y-y2,z-x;
poly g   = phi(f); 
g;
==> -x2y+yz2+x2y2-y2z2+x5+5x4z+10x3z2+10x2z3+5xz4+z5+y7-7y8+21y9-35y10
==> -x11+35y11+11x10z-55x9z2+165x8z3-330x7z4+462x6z5-462x5z6+330x4z7
==> -165x3z8+55x2z9-11xz10+z11-21y12+7y13-y14
classify(g);
==> The singularity ... is R-equivalent to T[p,q,r]=T[5,7,11]
\end{verbatim}

Ingredients for the classification of singularities:

\begin{enumerate}
\item standard bases for local and global orderings;
\item computation of invariants (Milnor number, determinacy, $\ldots$);
\item generalised Morse lemma;
\item syzygies for local orderings.
\end{enumerate}

Beyond classification by normal forms, the construction of moduli spaces for
singularities, for varieties or for vector bundles is a pretentious goal,
theoretically as well as computational.  First steps towards this goal for
singularities have been undertaken in \ccite{Ba} and \ccite{FrK}.
\vspace{0.5cm}

\subsection{deformations}

Consider a singularity $(V,0)$ given by power series $f_1,\dots, f_k \in
K\langle x_1, \dots, x_n\rangle$.  The idea of deformation theory is to
perturb the defining functions, that is to consider power series $F_1(t,x),
\dots,$ $ F_k(t,x)$ with $F_i(0,x) = f_i(x)$, where $t\in S$ may be considered
as a small parameter of a parameter space $S$ (containing 0).

For $t \in S$ the power series $f_{i,t}(x) = F_i(t,x)$ define a singularity
$V_t$, which is a perturbation of $V = V_0$ for $t\not= 0$ close to $0$.  It
may be hoped that $V_t$ is simpler than $V_0$ but still contains enough
information about $V_0$.  For this hope to be fulfilled, it is, however,
necessary to restrict the possible perturbations of the equations to
\textit{flat}\/ perturbations, which are called deformations.

Grothendieck's criterion of flatness states that the perturbation given by the
$F_i$ is \textbf{flat}\/ if and only if any relation between the $f_i$, say
\[
\sum r_i(x) f_i(x) = 0,
\]
lifts to a relation   
\[
\sum R_i(t,x) F_i(t,x)=0,
\]
with $R_i(x,0) = r_i(x)$. 
Equivalently, for any generator
$(r_1, \dots, r_k)$ of $\syz(f_1, \dots, f_k)$ there exists an
element $(R_1, \dots, R_k) \in \syz(F_1, \dots, F_k)$
satisfying $R_i(0,x) = r_i(x)$.   Hence, syzygies with respect to local
orderings come into play.   

There exists the notion of a semi--universal
deformation of $(V,0)$ which contains essentially all information about all
deformations of $(V,0)$.

For an isolated hypersurface singularity $f(x_1, \dots, x_n)$ the
semi--universal deformation is given by
\[
F(t,x) = f(x) + \sum^\tau_{j=1} t_j g_j(x),
\] 
where $1 =: g_1, g_2, \dots, g_\tau$ represent a $K$--basis of the
Tjurina algebra 
\[
K\langle x \rangle /\langle f,  \partial f/\partial x_1, \dots,
\partial f/\partial x_n\rangle,
\]
$\tau = \dim_K K\langle x \rangle/\langle f, \partial f/\partial x_1, \dots,
\partial f/\partial x_n\rangle$
being the Tjurina number.

To compute $g_1,\ldots,g_\tau$ we only need to compute a standard basis of the
ideal $\langle f, \tfrac{\partial f}{\partial x_1}, \dots, \tfrac{\partial
  f}{\partial x_n}\rangle$ with respect to a local ordering and then compute a
basis of $K[x]$ modulo the leading monomials of the standard basis.  For
complete intersections we have similar formulas.

\begin{center}
  \unitlength1cm
  \begin{picture}(13,6)
\put(0,1.2){\includegraphics[clip,height=5cm,width=5cm]{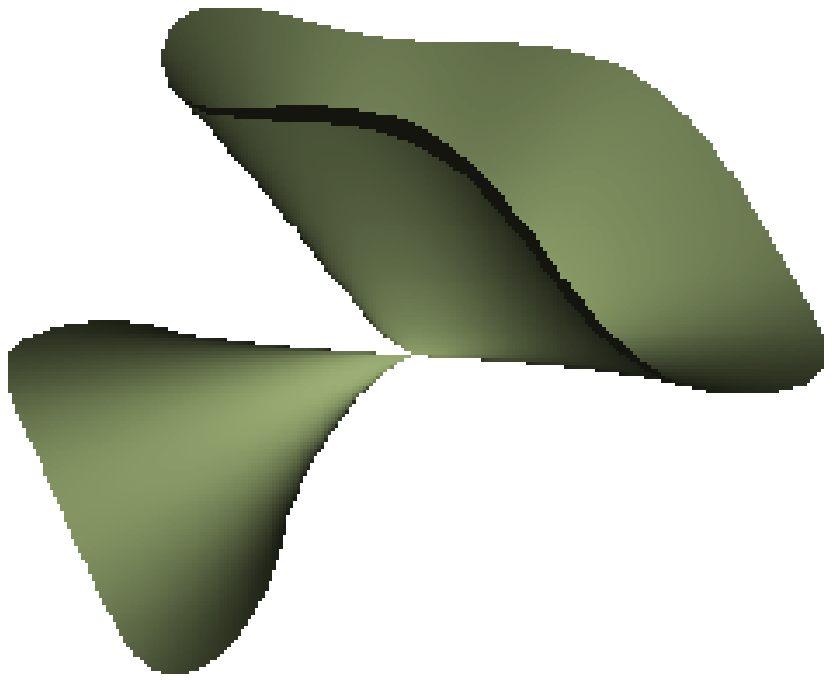}}
\put(8,1){\includegraphics[clip,height=5cm,width=5cm]{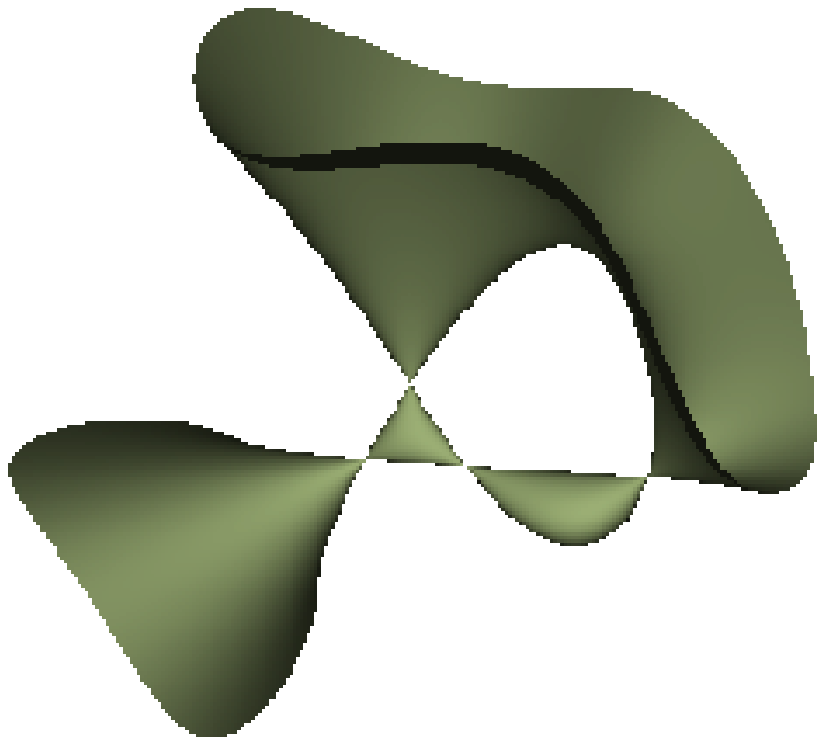}}
\put(4,0.5){Deformation of $E_7$ in $4A_1$}
  \end{picture}
\end{center}

For non--hypersurface singularities, the semi--universal deformation is much
more complicated and up to now no finite algorithm is known in general.
However, there exists an algorithm to compute this deformation up to arbitrary
high order cf.\ \ccite{Ll,Ma1}, which is implemented in SINGULAR.

As an example we calculate the base space of the semi--universal deformation
of the normal surface singularity, being the cone over the rational normal
curve $C$ of degree 4, parametrised by $t \mapsto (t,t^2,t^3,t^4)$.

Homogeneous equations for the cone over $C$ are given by the $2 \times
2$--minors of the matrix:
\[
m = \binom{x\;\;y\;\;z\;\;u}{y\;\;z\;\;u\;\;v} \in \Mat_{2\times
  4}(K[x,y,z,u,v]).
\]

\noindent SINGULAR commands for computing the semi--universal deformation:

\begin{minipage}[c][6cm]{7cm}
 \begin{verbatim}
LIB "deform.lib";
ring r         = 0,(x,y,z,u,v),ds;
matrix m[2][4] = x,y,z,u,y,z,u,v;
ideal f        = minor(m,2);
versal(f);                         
setring Px;
Fs;
==> Fs[1,1]=-u2+zv+Bu+Dv
==> Fs[1,2]=-zu+yv-Au+Du
==> Fs[1,3]=-yu+xv+Cu+Dz
==> Fs[1,4]=z2-yu+Az+By
==> Fs[1,5]=yz-xu+Bx-Cz
==> Fs[1,6]=-y2+xz+Ax+Cy
Js;
==> Js[1,1]=BD
==> Js[1,2]=AD-D2
==> Js[1,3]=-CD  
\end{verbatim}
\end{minipage}
\begin{minipage}[c][6cm]{7cm}
  \unitlength1cm
  \begin{picture}(5,7)
  \put(-0.5,-1.5){\includegraphics[clip,width=7.5cm]{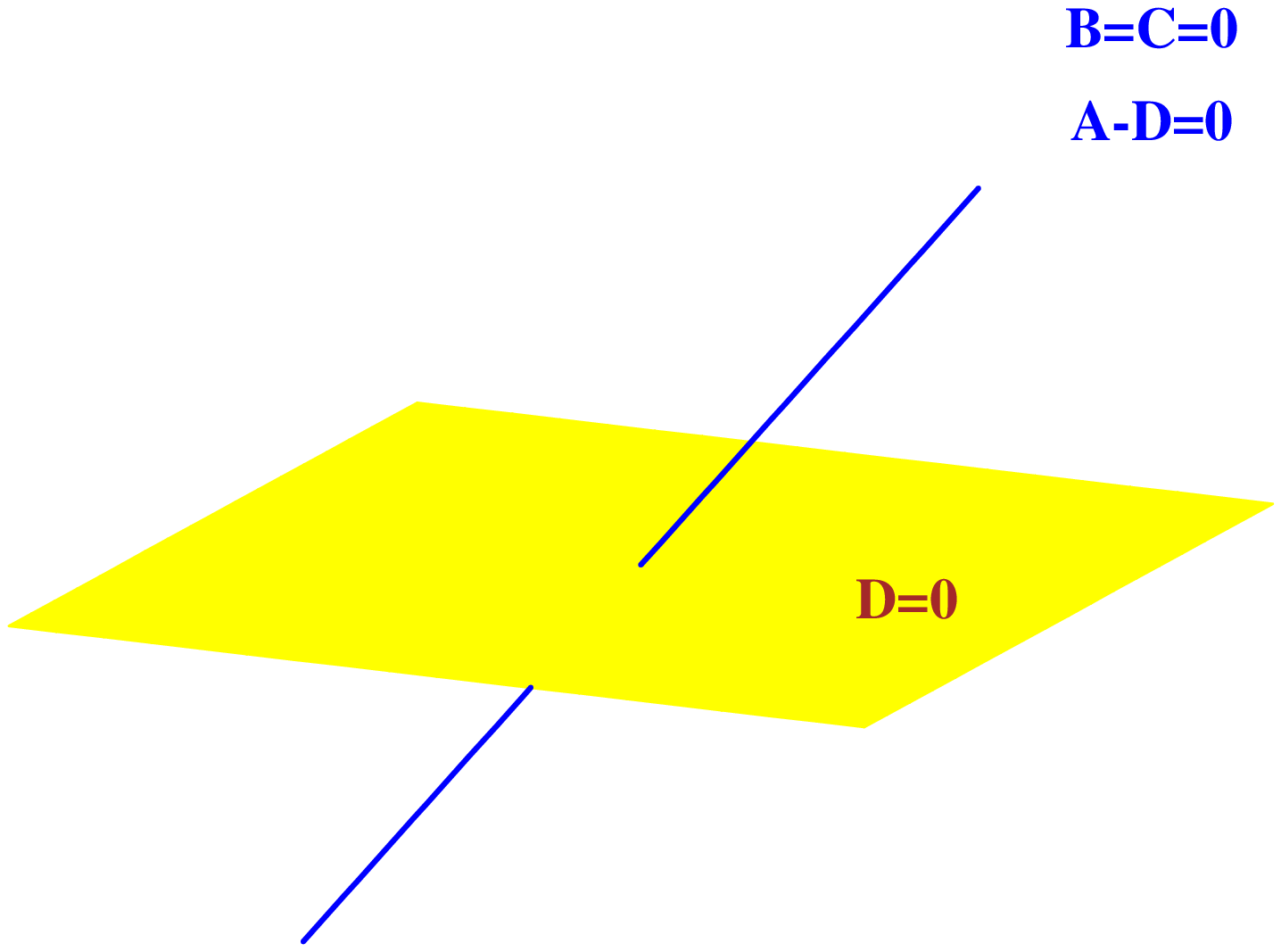}}
  \end{picture}
\end{minipage}
\medskip

The ideal $Js =\langle BD, AD-D^2,-CD\rangle \subset K[A,B,C,D]$ defines the
required base space which consists of a 3--dimensional component $(D=0)$ and a
transversal 1--dimensional component $(B=C=A-D=0)$.   This was the first
example, found by Pinkham, of a base space of a normal surface having several
components of different dimensions.

The full versal deformation is given by the map (\verb?Fs? and \verb?Js? as
above)
\[
K[[A,B,C,D]]/J_s \lra K[[A,B,C,D,x,y,z,u,v]]/J_s + F_s.
\]

Although, in general, the equations for the versal deformation are formal
power series, in many cases of interest (as in the example above) the
algorithm terminates and the resulting ideals are polynomial.  \medskip

Ingredients for the semi--universal deformation algorithm:

\begin{enumerate}
\item Computation of standard bases, normal forms and resolutions for local
  orderings;

\item computation of Ext groups (cf.\ 4.1) for computing infinitesimal
  deformations and obstructions;

\item computation of Massey products for determining obstructions to
  lift, recursively, infinitesimal deformations of a given order to higher
  order; 

\item one of the main difficulties in point 3 is the necessity to compute a
  completely reduced normal form with respect to a local ordering.   In
  general, such a normal form exists only as formal power series.   In the
  present situation, however, the reduction has to be carried out only for a
  subset of the variables  in a fixed degree
  and, hence, the complete reduction is finite.
\end{enumerate}

\subsection{the monodromy}

Let $f \in \C\{x_1, \dots, x_n\}$ be a convergent power series (in practice a
polynomial) with isolated singularity at 0 and $\mu = \dim_\C \C\{x\}/\langle
f_{x_1}, \dots, f_{x_n}\rangle$ the Milnor number of $f$.

Then $f$ defines in an $\varepsilon$-ball $B_\varepsilon$ around 0 a
holomorphic function to $\C$, $f : B_\varepsilon \lra \C$.

The simple, counterclockwise path $\gamma$ in $\C$ around $0$ induces a
$C^\infty$--diffeomorphism of $X_t\; (t \not= 0)$ (as indicated in the 
figure) and an automorphism of the singular cohomology group $H^n(X_t,
\C)$ which is, by a theorem of Milnor, a $\mu$--dimensional $\C$--vector
space.  This automorphism
\[
T : H^n(X_t, \C) \stackrel{\cong}{\lra} H^n(X_t, \C)
\]
is called the local \textbf{Picard--Lefschetz monodromy} of $f$.   We address
the problem of computing the Jordan normal form of $T$.

\begin{center}
  \vspace*{0.3cm}
   \unitlength1cm
     \begin{texdraw}
       \drawdim{cm} 
       \move (5 7) \lcir r:2 
       \move (5 2) \lellip rx:2 ry:1 \fcir f:0 r:0.05
                   \linewd 0.03 \lpatt(0.1 0.14) \lellip rx:1 ry:0.5
                   \lpatt()
       \htext (5.2 1.8){$0$}
       \move (5.8 1.75) \arrowheadtype t:V \avec (5.95 1.9)
       \move (6 2) \fcir f:0 r:0.05       
       \htext (6.2 1.8){$t$}
       \htext (5.8 2.4){$\gamma$ \hspace{0.8cm}= path around $0$}
       \htext (5.2 4){$f$}
       \htext (5 9.1){$X_0=B_\varepsilon\cap f^{-1}(0)$}
       \htext (6.1 4.8){$X_t=B_\varepsilon\cap f^{-1}(t)$}
       \move (5 9) \linewd 0.03 \clvec (5.5 8.5)(5.1 7)(4.5 7)
                                 \clvec (5.4 7)(5.3 5.8)(5 5)
       \move (6 8.73) \clvec (7 7.5)(5.5 6.3)(6 5.27)
       \move (5 4.7) \arrowheadtype t:F \avec (5 3.3)
       \move (6.3 7.7) \arrowheadtype t:V \lpatt(0.1 0.14) \cavec (3.7 
8.2)(3.7 6.8)(6.2 7.2)
       \move (5.9 6.1) \arrowheadtype t:V \lpatt(0.1 0.15) \cavec (3.6 
6.7)(3.5 5.2)(5.9 5.9)
     \end{texdraw}
\end{center}

The first important theorem is the

\noindent \textbf{Monodromy theorem} (Deligne 1970, Brieskorn 1971): \textit{The 
  eigenvalues of $T$ are roots of unity, that is, we have
\[
T = e^{2\pi i M},
\]
where $M$ is a complex matrix with eigenvalues in $\Q$.}

Hence, we are left with the problem of computing the Jordan normal form of
$M$.
\medskip

It is not at all clear that the purely topological definition of $T$ allows an
algebraic and computable interpretation.   The first hint in this direction is
that we can compute $\dim_\C H^n(X_t, \C)$, according to Milnor's theorem,
algebraically by the formula for $\mu$ given above.

Since $X_t$ is a complex Stein manifold, its complex cohomology
can be computed, via the holomorphic de Rham theorem, with holomorphic
differential forms, which is the starting point for computing the monodromy.

To cut a long story short, we just mention, cf.\ \ccite{Br,Gre1}
that \vspace{-0.5cm}

\begin{align*}
  H' & = \Omega^n/df  \wedge \Omega^{n-1} + d \Omega^{n-1},\\
H^{\prime\prime} & = \Omega^{n+1}/df  \wedge d \Omega^{n-1}
\end{align*}
are free $\C\{t\}$--modules (via $f^\ast : \C\{t\} \lra \C\{f\} \subset
\C\{x\}$) of rank equal to $\mu$.   Here $(\Omega^\bullet,d)$ denotes the
complex of holomorphic differential forms in $(\C^n,0)$.   $H^\prime$ and
$H^{\prime\prime}$ are called Brieskorn lattices.

We define the local \textbf{Gau{\ss}--Manin connection} of $f$ as
\begin{gather*}
  \bigtriangledown : df  \wedge H^\prime =  df \wedge
  \Omega^n/df  \wedge d \Omega^{n-1} \lra H^{\prime\prime},\\
\bigtriangledown [df \wedge \omega] = [d \omega].
\end{gather*}

Note that $\bigtriangledown(df \wedge H^\prime) \not\subset df \wedge
H^\prime$, that is, $\bigtriangledown$ has a pole at $0$.  Tensoring with
$\C\,(t)$, the quotient field of $\C\{t\}$, we can extend $\bigtriangledown$
to a meromorphic connection
\[
\bigtriangledown : H^{\prime\prime} \underset{\C\{t\}}{\otimes} \C(t) \lra
H^{\prime\prime} \underset{\C\{t\}}{\otimes} \C(t)
\]
(since df $\wedge H^\prime \otimes \C(t) = H^{\prime\prime}\otimes \C(t)$)
using the Leibnitz rule $\bigtriangledown (\omega y) = \bigtriangledown
(\omega) y + \omega dy/dt$.

With respect to a basis $\omega_1, \dots, \omega_\mu$ of $H^{\prime\prime}$ we
have $\bigtriangledown (\omega_i) = \underset{j}{\sum} a_{j i} \omega_j$ and,
for any $\omega = \underset{i}{\sum} \omega_i y_i$, $\bigtriangledown (\omega)
= \underset{i,j}{\sum} a_{ji} y_i + \underset{i}{\sum} \omega_i dy_i/dt$.
Hence, the kernel of $\bigtriangledown$, together with a basis of
$H^{\prime\prime}$, is the same as the solutions of the system of rank $\mu$
of ordinary differential equations
\[
\dfrac{dy}{dt} = -Ay,\qquad A = (a_{ij}) \in \text{ Mat}\bigl(\mu \times
\mu, \C(t)\bigr)
\]
in a neighbourhood of $0$ in $\C$.  The connection matrix, $A = \underset{i
  \ge -p}{\sum} A_i t^i$, $A_i \in \Mat(\mu \times \mu, \C)$, has a pole at $t
= 0$ and is holomorphic for $t \not= 0$.  If $\phi_t = (\phi_1, \dots,
\phi_\mu)$ is a fundamental system of solutions at a point $t \not= 0$, then
the analytic continuation of $\phi_t$ along the path $\gamma$ transforms
$\phi_t$ into another fundamental system $\phi'_t$ which satisfies $\phi'_t =
T_\bigtriangledown \phi_t$ for some matrix $T_\bigtriangledown \in \text{
  GL}(\mu, \C)$.
\smallskip

\noindent \textbf{Fundamental fact} (Brieskorn, 1971): \textit{The
  Picard--Lefschetz   monodromy $T$ coincides with the monodromy
  $T_\bigtriangledown$ of the Gau{\ss}--Manin connection.}
\medskip

Brieskorn used this fact to describe in \ccite{Br} the essential steps for an
algorithm to compute the characteristic polynomial of $T$.  Results of Gerard
and Levelt allowed the extension of this algorithm to compute the Jordan
normal form of $T$, cf.\ \ccite{GL}.  An early implementation by Nacken in
Maple was not very efficient.  Recently, \ccite{Schu} implemented an
improved version in SINGULAR which is able to compute interesting examples.

The algorithm uses another basic theorem, the
\vspace{0.5cm}

\noindent \textbf{Regularity Theorem} (Brieskorn, 1971): \textit{The Gau{\ss}--Manin
  connection has a regular singular point at $0$, that is, there exists a
  basis of some lattice in $H^{\prime\prime} \otimes \C(t)$ such that the
  connection matrix $A$ has pole of order 1.}
\medskip

Basically, if $A = A_{-1} t^{-1} + A_0 + A_1 t + \cdots$ has a simple pole,
then $T = e^{2\pi i} A_{-1}$ is the monodromy (this holds if the eigenvalues of
$A_{-1}$ do not differ by integers which can be achieved algorithmically).
\smallskip

\begin{samepage}
\noindent SINGULAR example:
\begin{verbatim}
LIB "mondromy.lib";
ring R   = 0,(x,y),ds;
poly f   = x2y2+x6+y6;   //example of A'Campo (monodromy is not diagonalisable)
matrix M = monodromy(f);
print(jordanform(M));
==> 1/2,1,  0,  0,  0,  0,  0,0,0,0,  0,  0,  0, 
==> 0,  1/2,0,  0,  0,  0,  0,0,0,0,  0,  0,  0, 
==> 0,  0,  2/3,0,  0,  0,  0,0,0,0,  0,  0,  0, 
==> 0,  0,  0,  2/3,0,  0,  0,0,0,0,  0,  0,  0, 
==> 0,  0,  0,  0,  5/6,0,  0,0,0,0,  0,  0,  0, 
==> 0,  0,  0,  0,  0,  5/6,0,0,0,0,  0,  0,  0, 
==> 0,  0,  0,  0,  0,  0,  1,0,0,0,  0,  0,  0, 
==> 0,  0,  0,  0,  0,  0,  0,1,0,0,  0,  0,  0, 
==> 0,  0,  0,  0,  0,  0,  0,0,1,0,  0,  0,  0, 
==> 0,  0,  0,  0,  0,  0,  0,0,0,7/6,0,  0,  0, 
==> 0,  0,  0,  0,  0,  0,  0,0,0,0,  7/6,0,  0, 
==> 0,  0,  0,  0,  0,  0,  0,0,0,0,  0,  4/3,0, 
==> 0,  0,  0,  0,  0,  0,  0,0,0,0,  0,  0,  4/3
\end{verbatim}
\end{samepage}
Ingredients for the monodromy algorithm:

\begin{enumerate}
\item Computation of standard bases and normal forms for local orderings;
\item computation of Milnor number;
\item Taylor expansion of units in $K[x]_{\langle x \rangle}$ up to
  sufficiently high order;
\item computation of the connection matrix on increasing lattices in
  $H^{\prime\prime} \otimes \C(t)$ up to sufficiently high order (until
  saturation) by linear algebra over $\Q$\,;
\item computation of the transformation matrix to a simple pole by linear
  algebra over $\Q$\,;
\item factorisation of univariate polynomials (for Jordan normal form). 
\end{enumerate}

The most expensive parts are certain normal form computations for a local
ordering and the linear algebra part because here one has to deal iteratively
with matrices with several thousand rows and columns.  It turned out that the
SINGULAR implementation of modules (considered as sparse matrices) and the
Buchberger inter-reduction is sufficiently efficient (though not the best
possible) for such tasks.  

\section{Computer algebra solutions to singularity problems}

We present three examples which demonstrate, in a somewhat typical way, the use
of computer algebra as stated in the preface:

\begin{enumerate}
\item producing counter examples;
\item providing evidence and prompting proofs for new theorems;
\item constructing interesting explicit examples;
\end{enumerate}

\subsection{exactness of the Poincar\'{e} complex}

The first application is a counterexample to a conjectured generalisation of a
theorem of \ccite{Sa} which says that, for an isolated hypersurface
singularity, the exactness of the Poincar\'{e} complex implies that the
defining polynomial is, after some analytic coordinate change, weighted
homogeneous.  \smallskip

\noindent \textbf{Theorem} \cite{Sa}:

\textit{
If $f : \C^{n+1} \lra \C$ has an isolated singularity at $0$, then the
following are equivalent:}

\begin{enumerate}
\item \textit{$X = f^{-1} (0)$ is weighted homogeneous for a suitable choice of
  coordinates.}
  
\item \textit{$\mu = \tau$ where $\mu = \dim_\C \C\{x\}/\left(\dfrac{\partial
      f}{\partial x_i}\right)$ is the Milnor number and\\
 $\tau = \dim_\C
  \C\{x\}/\left(f,\dfrac{\partial f}{\partial x_i}\right)$ the Tjurina number.}

\item \textit{The holomorphic Poincar\'{e} complex}
\[
0 \lra \C \lra \ko_X \overset{d}{\lra} \Omega^1_X \overset{d}{\lra} \Omega^2_X
\lra \dots \lra \Omega^n_X \lra 0
\]
\textit{is exact.}
\end{enumerate}

A natural problem is whether the theorem holds also for complete intersections
$X = f^{-1} (0)$ with $f = (f_1, \dots, f_k) : \C^{n+k} \lra \C^k$.  Again we
have a Milnor number $\mu$ and a Tjurina number $\tau$,
\begin{align*}
  \mu & = \sum^k_{i=1} (-1)^{i-1} \dim_\C \C \{x\}/\left(f_1,
    \dots, 
  f_{i-1}, \Biggl\lvert \dfrac{\partial(f_1, \dots, f_i)}{\partial(x_{j_1},
    \dots, x_{j_i})}\Biggr\lvert \right)\\
\tau & = \dim_\C \C\{x\}^k/\bigl(f_1, \dots, f_k) \C\{\underline{x}\}^k
+ Df (\C\{\underline{x}\}^{n+k})\bigr).
\end{align*}
\smallskip

\noindent {\bf Theoretical reduction} \cite{GMP}:

\textit{If $X$ is a complete intersection of dimension 1, then (1)
  $\Leftrightarrow$ (2) $\Rightarrow$ (3).}  \\
\textit{\noindent If $k = 2$, then (3)
  $\Rightarrow$ (2) if $\mu = \dim_\C \Omega^2_X - \dim_\C \Omega^3_X$ and if
  $f_1, f_2$ are weighted homogeneous.}
\medskip

\ccite{PS} showed that (3) $\Rightarrow$ (2) does {\bf not} hold in general:
\[
f_1 = xy + z^{\ell-1},\;\; f_2 = xz + y^{k-1} + yz^2 \quad (4 \le \ell \le
k,\; k \ge 5)
\]
is a counterexample.

The proof uses an implementation of the standard basis algorithm in a
forerunner of SINGULAR and goes as follows:

\begin{enumerate}
\item Compute $\mu, \dim_\C \Omega^2_X, \dim_\C \Omega^3_X$ to show that
  $\Omega^\bullet_X$ is exact;

\item compute $\tau$.\\
One obtains $\mu = \tau + 1$, that is, $X$ is not weighted homogeneous.
\end{enumerate}

To do this we must be able to compute standard bases of modules over local
rings.

The counterexample was found through a computer search in a list of
singularities classified by \ccite{Wa}.

\subsection{Zariski's multiplicity conjecture}

The attempt to find a counterexample to Zariski's multiplicity conjecture ---
which says that the multiplicity (lowest degree) of a power series is an
invariant of the embedded topological type --- led, finally after many
experiments and computations, to a partial proof of this conjecture. For this,
an extremely fast standard basis computation for 0--dimensional ideals in a
local ring was necessary.

The following question was posed by \ccite{Za} in his retiring address
to the AMS in 1971.

Let $f = \sum c_\alpha x^\alpha \in \C\{x_1, \dots, x_n\},\; f(0) = 0$, be a
hypersurface singularity, and let $\mult(f) := \min\{|\alpha|\; \Bigl\lvert\;
c_\alpha \not= 0\Bigr.\}$ be the multiplicity.

We say that $f$ and $g$ are topological equivalent, $f
\overset{\text{top}}{\sim} g$, if there is a homeomorphism 
\[
(B, f^{-1}(0) \cap
B,0) \overset{\sim}{\lra} (B, g^{-1}(0) \cap B,0)
\]

\begin{center}
\unitlength1cm
\begin{picture}(12,4.5)
\put(0,0){\includegraphics[clip,angle=90,height=5cm,width=5cm]{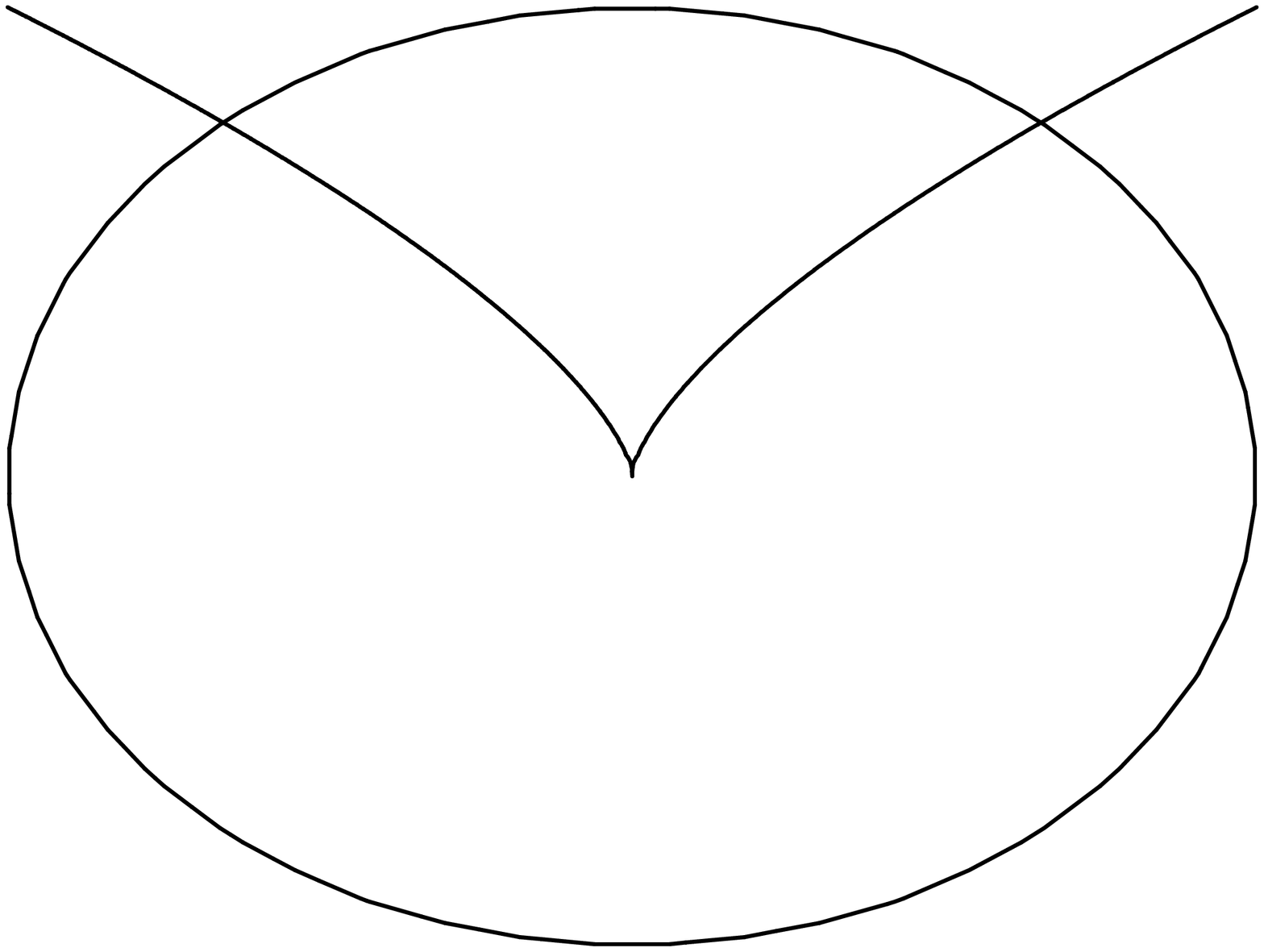}}
\put(6,0){\includegraphics[clip,height=5cm,width=5cm]{Bilder/kuspe-mit-Milnorkreis.ps}}
\put(-0.5,4){$f^{-1}(0)$}
\put(2.5,2){0}
\put(5.3,2){$\overset{\sim}{\lra}$}
\put(8.3,2){0}
\put(10.5,4){$g^{-1}(0)$}
\end{picture}
\end{center}

\noindent {\bf Zariski's conjecture} may be stated as:
$f\overset{\text{top}}{\sim} g \Rightarrow \mult(f) = \mult (g)$.
\medskip
 
 The result is known to be true for curves (Zariski, L\^{e}) and weighted
 homogeneous singularities \ccite{Gre3,OS}.

Our attempt to find a counterexample was as follows:\\
Consider deformations of $f = f_0$:
\[
f_t(x) = f(x) + tg(x,t), \quad |t| \text{ small}.
\]
Then use the theoretical fact  proved by L\^{e} and Ramanujam:
\[
f_0\overset{\text{top}}{\sim} f_t \Rightarrow \mu(f_0) = \mu(f_t)
\]
(``$\Leftarrow$'' holds also, except for $n = 3$, where the answer is still
unknown) where $\mu(f_0)$ respectively $\mu(f_t)$ are the
Milnor numbers. 

We tried to construct a deformation $f_t$ of $f_0$ where the multiplicity $\mult(f_t)$
drops but the Milnor number $\mu(f_t)$ is constant.

Our candidates $(a, b, c \in \N)$ came from a heuristical investigation of the
Newton diagram, one being the following series:
 \[
f_t = x^a + y^b + z^{3c} + x^{c+2} y^{c-1} + x^{c-1} y^{c-1} z^3 + x^{c-2} y^c
(y^2 + tx)^2,\; a,b,c \in \N.
\]
Obviously, the multiplicity drops.   Computing $\mu$ with SINGULAR, we obtain
for $(a, b, c) = (37, 27, 6)$: $\mu(f_0) = 4840$, $\mu(f_t) = 4834$, thus
$f_0$ and $f_t$ are (unfortunately) not topologically equivalent.

Since the Milnor numbers of possible counter examples have to  be very big, we
need an extremely efficient implementation of standard bases.   For this, the
``highest corner'' method of \ccite{GP1} was essential.

Trying many other classes of examples, we did not succeed in finding a counter
example.   However, an analysis of the examples led to the following
\smallskip

\noindent {\bf Partial proof of Zariski's conjecture} \cite{GP1}:

\noindent \textit{Zariski's conjecture is true for deformations of the form}
\[
f_t = g_t(x,y) + z^2 h_t(x,y),\quad \mult(g_t) < \mult(f_0).
\]
There is also an invariant characterisation of the deformations of the above
kind.   The general conjecture is, up to today, still open.

\subsection{curves with maximal number of singularities}

Let $C \subset \P^2_\C$ be an irreducible projective curve of degree $d$ and
$f(x,y) = 0$ a local equation for the germ $(C,z)$.   Let $\mu(C,z) =
\dim_\C \C\{x,y\}/(f_x,f_y)$ be the Milnor number of $C$ at $z$.

Since the genus of $C$, $g(C) = \tfrac{(d-1)(d-2)}{2} - \delta(C)$ is
non--negative (where $\delta(C) = \underset{z\in C}{\sum} \delta(C,z)$,
$\delta(C,z) = \dim_\C \bar{R}/R$, $R = \C\{x,y\}/\langle f\rangle$ and
$\bar{R}$ the normalisation of $R$), $C$ can have, at most, $(d-1)(d-2)/2$
singularities.

It is a classical and interesting problem, which is still in the centre of
theoretical research, to study the variety $V = V_d(S_1, \dots, S_r)$ of
(irreducible) curves $C \subset \P^2_\C$ of degree $d$ having exactly $r$
singularities of prescribed (topological or analytical) type $S_1, \dots,
S_r$.  Among the most important questions are:
\begin{itemize}
\item Is $V \not= \emptyset$  (existence problem)?
\item Is $V$ irreducible (irreducibility problem)?
\item Is $V$ smooth of expected dimension ($T$--smoothness problem)?
\end{itemize}
A complete answer is only known for nodal curves, that is, for $V_d(r) =
V_d(S_1, \dots, S_r)$ with $S_i$ ordinary nodes ($A_1$--singularities):

Severi (1921):  $V_d(r) \not= \emptyset$ and $T$--smooth $\Leftrightarrow r
\le \tfrac{(d-1)(d-2)}{2}$.

Harris (1985): $V_d(r)$ is irreducible (if $\not= \emptyset$).

Even for cuspidal curves a sufficient and necessary answer to any of the above
questions is unknown.

\unitlength1cm
\begin{picture}(13,6.5)
    \put(0,0){
       \begin{minipage}[b]{6cm}
         \begin{center}
            {\includegraphics[trim=20 20 20 20, clip]{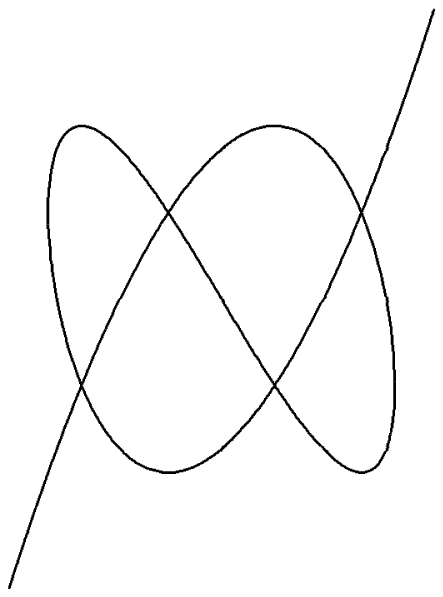}}
         \end{center}
A $4$-nodal plane curve of degree 5, with equation
$x^5-\frac{5}{4}x^3+\frac{5}{16}x-\frac{1}{4}y^3+\frac{3}{16}y=0$,
which is a deformation of $E_8\; : \; x^5-y^3=0$.
       \end{minipage}
       }
\put(7,0){
\begin{minipage}[b]{6cm}
\begin{center}
    \includegraphics[clip,height=5cm,width=5cm]{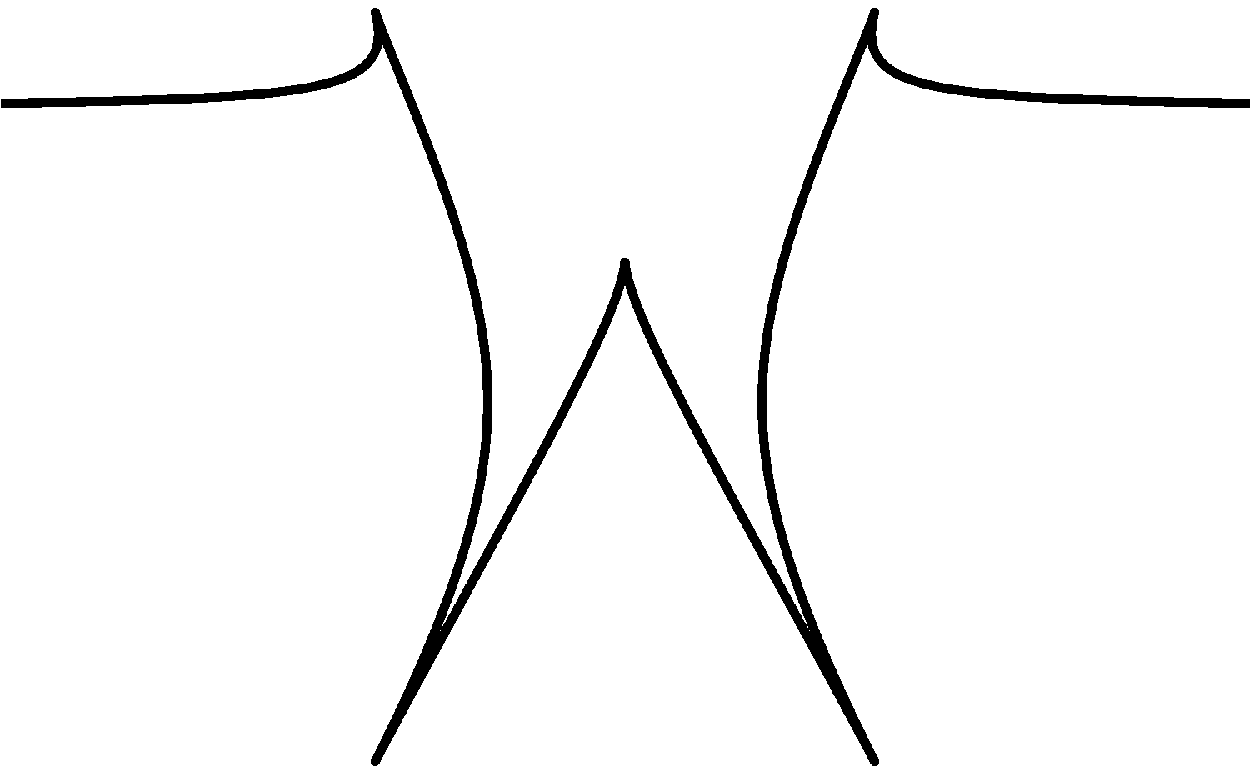}
\end{center}
A plane curve of degree $5$ with 5 cusps, the maximal possible number.   It
has the equation $\frac{129}{8}x^4y-\frac{85}{8}x^2y^3+\frac{57}{32}y^5-20x^4-\frac{21}{4}x^2y^2+\frac{33}{8}y^4-12x^2y+\frac{73}{8}y^3+32x^2=0$.
\end{minipage}}
\end{picture}
\vspace{1cm}

Concerning arbitrary (topological types of) singularities, we have the
following existence theorem, which is, with respect to the exponent
of $d$, asymptotically optimal.
\vspace{0.5cm}

\noindent \textbf{Theorem:} \cite{GLS,Lo}.

\noindent {\it $V_d(S_1, \dots, S_r) \not=\emptyset$ if $\sum^r_{i=1} \mu(S_i)
  \le \dfrac{(d+2)^2}{46}$ and two additional conditions for the five
  ``worst'' singularities.

In case of only one singularity we have the slightly better sufficient
condition for existence,}
\[
\mu(S_1) \le \dfrac{(d-5)^2}{29}.
\]

The theorem is just an existence statement, the proof gives no hint how to
produce any equation.  Having a method for constructing curves of low degree
with many singularities, Lossen was able to produce explicit equations.  In
order to check his construction and improve the results, he made extensive use
of SINGULAR to compute standard bases for global as well as for local
orderings.  One of his examples is the following: \smallskip

\noindent \textbf{Example}: \cite{Lo}\quad The irreducible curve with affine
equation $f(x,y) = 0$, 
\begin{align*}
  f(x,y) = y^2 - 2y(x^{10} & + \dfrac{1}{2} x^{9} y^2 - \dfrac{1}{8} x^{8}
  y^4 + \dfrac{1}{16} x^{7} y^6 - \dfrac{5}{128} x^{6} y^8 + \dfrac{7}{256}
  x^5
  y^{10}\\
  & - \dfrac{21}{1024} x^4 y^{12} + \dfrac{33}{2048} x^3 y^{14} -
  \dfrac{429}{32768} x^2 y^{16} + \dfrac{715}{65536} x y^{18}\\
  & - \dfrac{2431}{262144}  y^{20})  + x^{20} + x^{19}y^2
\end{align*}
has degree 21 and an $A_{228}$--singularity $(x^2-y^{229} = 0)$ as its only
singularity.

In order to verify this, one may proceed, using SINGULAR, as follows:
\begin{verbatim}
ring s = 0,(x,y),ds;
poly f = y2-2x10y-x9y3+1/4x8y5-1/8x7y7+5/64x6y9-7/128x5y11+21/512x4y13 
        -33/1024x3y15+429/16384x2y17+x20-715/32768xy19+x19y2+2431/131072y21;
matrix Hess = jacob(jacob(f));         //the Hessian matrix of f
print(subst(subst(Hess,x,0),y,0));     //the Hessian matrix for x=y=0
==> 0,0,
==> 0,2 
vdim(std(jacob(f)));                  //the Milnor number of f
==> 228
\end{verbatim}

Since the rank of the Hessian at 0 is 1, $f$ has an $A_k$ singularity at 0; it
is an $A_{228}$ singularity since the Milnor number is 228.  To show that the
projective curve $C$ defined by $f$ has no other singularities, we have to
show that $C$ has no further singularities in the affine part and no
singularity at infinity. The second assertion is easy, the first follows from
\[
        dim_\C(K[x,y]_{\langle x,y\rangle}/\langle\text{jacob}(f),f\rangle =
        \dim_\C(K[x,y]/\langle \text{jacob}(f),f\rangle,
\]
confirmed by SINGULAR:
\begin{verbatim}
vdim(std(jacob(f)+f));
==> 228                     //multiplicity of Sing(C) at 0 (local ordering)
ring r = 0,(x,y),dp;
poly f = fetch(s,f);
vdim(std(jacob(f)+f));
==> 228                     //multiplicity of Sing(C)      (global ordering)
\end{verbatim}

\section{What else is needed}

In this survey I could only touch on a few topics where computer algebra has
contributed to mathematical research.  Many others have not been mentioned,
although there exist powerful algorithms and efficient implementations.  In
the first place, the computation of invariant rings for group actions of
finite \cite{Stu,Ke,DJ}, reductive \cite{Der} or some uni--potent \cite{GPS}
groups belong here.  Computation of invariants have important applications for
explicit construction of moduli spaces, for example, for vector bundles or for
singularities \cite{FrK,Ba} but also for dynamical systems with symmetries
\cite{Ga}.  Libraries for computing invariants are available in SINGULAR.  Available is also the Puiseux expansion (even better,
the Hamburger--Noether expansion, cf.\ \ccite{Lam} for description of an
implementation) of plane curve singularities.  The latter is one of the few
examples of an algorithm in algebraic geometry where Gr\"obner bases are not
needed.

The applications of computer algebra and, in particular, of Gr\"obner bases in
projective algebraic geometry are so numerous that I  can only refer to the
textbooks of Cox--Little--O'Shea, Eisenbud and Vasconcelos and the literature
cited there.   The applications include classification of varieties and vector
bundles, cohomology, moduli spaces and fascinating problems in enumerative
geometry. 

However, there are also some important problems for which an algorithm is
either not known or not yet implemented (for further open problems see also \ccite{Ei1}):

\begin{enumerate}
\item Resolution of singularities

This is one of the most important tools for treating singular varieties.   At
least three approaches seem to be possible:

 For surfaces we have Zariski's method of successive normalisation
  and blowing up points and the Hirzebruch--Jung method of resolving the
  discriminant curve of a projection.

 For arbitrary varieties, new methods of Bierstone, Milman and
  Villamajor provide a constructive approach to resolution in the spirit of
  Hironaka.   First attempts in this direction have been made by Schicho.   

\item Computation in power series rings

This is a little vague since I do not mean to actually compute with infinite
power series, the input should be polynomials.   However, it would be highly
desirable to make effective use of the Weierstra{\ss} preparation theorem.
This is related to the problem of elimination in power series rings.

Moreover, no algorithm seems to be known to compute an algebraic
representative of the semi-universal deformation of an isolated singularity
(which is known to exist).   

Also, I do not know any algorithm for Hensel's lemma.

\item Dependence of parameters
  
  In this category falls, at least principally, the study of Gr\"obner bases
  over rings.  This has, of course, been studied, cf.\ \ccite{AL,Ka1}, 
  but I still consider the dependence of Gr\"obner bases on parameters as an
  unsolved problem (in the sense of an intrinsic or predictable description,
  if it exists).

In many cases, one is interested in finding equations for parameters
describing precisely the locus where certain invariants jump.   This is
related to the above problem since Gr\"obner bases usually only give a
sufficient but not necessary answer.

Mainly in practical applications of Gr\"obner bases to ``symbolic solving'',
parameters are real or complex numbers.   It would then be important to know,
for which range of the parameters the symbolic solution holds.

\item Symbolic--numeric algorithms
  
  The big success of numerical computations in real life problems seems to
  show that symbolic computation is of little use for such problems.  However,
  as is well--known, symbolic preprocessing of a system of polynomial (even
  ordinary and partial differential) equations may not only lead to much
  better conditions for the system to be solved numerically but even make
  numerical solving possible.
  
  There is continuous progress in this direction, cf.\ \ccite{CLO2,Moe,Ve},
  not only by Gr\"obner basis methods.  A completely different approach via
  multivariate resultants (cf.\ \ccite{CE}) has become favourable to several
  people due to the new sparse resultants by \ccite{Ge}.  However, an
  implementation in SINGULAR (cf.\ \ccite{We,Hil}) does not show superiority
  of resultant methods, at least for many variables against triangular set
  methods of either Lazard or M\"oller.  Nevertheless, much more has to be
  done.

The main disadvantage of symbolic methods in practical, real life applications
is its complexity.  Even if a system is able to return a symbolic answer in a
short time, this answer is often not humanly interpretable.  Therefore, a
symbolic simplification is necessary, either before, during, or after
generation.  Of course, the result must still be approximately correct.

This leads to the problem of validity of ``simplified'' symbolic computation.
A completely open subproblem is the validity resp.\ error estimation of
Gr\"obner basis computations with floating point coefficients.
 
The simplification problem means providing simple and humanly understandable
symbolic solutions which are approximately correct for numerical values in a
region which can be specified.  This problem belongs, in my opinion, perhaps to the
most important ones in connection with applications of computer algebra to
industrial and economical problems.

\item Non--commutative algorithms
  
  Before Gr\"obner bases were introduced by Buchberger, the so--called
  Ritt--Wu method, cf.\  \ccite{Ri}, was developed for symbolic computation in
  non--commutative rings of differential operators.  However, nowadays,
  commutative Gr\"obner bases are implemented in almost every major computer
  algebra system, whilst only few systems provide non--commutative algorithms.
  Standard bases for some non--commutative structures have been implemented in
  the system FELIX \cite{AK} as well as in an experimental version in SINGULAR;  the system Bergman and
  an extension called Anick can compute Gr\"obner bases and higher
  syzygies in the non--commutative case.
 
  Highly desirable are effective implementations for non--commutative
  Gr\"obner bases in the Weyl algebra, the Grassmannian, for D--modules or the
  enveloping algebra of a finite dimensional Lie algebra (the general theory
  being basically understood, cf.\ \ccite{Mo2,Uv,Ap}).
  
  The recent textbook of \ccite{SST} shows a
  wide variety of algorithms  for modules over the Weyl algebra and
  $D$--modules for which an efficient implementation is missing.
  
  But even classical algebraic geometry, as was shown, for example, by
  \ccite{KM}, has a natural embedding into
  non--commutative algebraic geometry.  A special case is known as
  quantisation, a kind of non--commutative deformation of a commutative
  algebra.

Providing algorithms and implementations for the use of computer algebra in
non--com\-mutative  algebraic geometry could become a task and challenge for a
new generation of computer algebra systems.
\end{enumerate}

\end{document}